\documentclass{article}
\usepackage{amssymb}

\usepackage{amsmath}

\newtheorem{theorem}{Theorem}

\newtheorem{corollary}[theorem]{Corollary}

\newtheorem{definition}[theorem]{Definition}
\newtheorem{example}[theorem]{Example}

\newtheorem{proposition}[theorem]{Proposition}
\newtheorem{remark}[theorem]{Remark}

\newenvironment{proof}[1][Proof]{\textbf{#1.} }{\ \rule{0.5em}{0.5em}}
\input{tcilatex}

\begin{document}

\title{Abelian equations and rank problems for planar webs }
\author{Vladislav V. Goldberg and Valentin V. Lychagin}
\maketitle

\begin{abstract}
We find an invariant characterization of planar webs of maximum rank. For $4$%
-webs, we prove that a planar $4$-web is of maximum rank three if and only
if it is linearizable and its curvature$\;$vanishes. This result leads to
the direct web-theoretical proof of the Poincar\'{e}'s theorem: a planar $4$%
-web of maximum rank is linearizable. We also find an invariant intrinsic
characterization of planar $4$-webs of rank two and one and prove that in
general such webs are not linearizable. This solves the Blaschke problem
``to find invariant conditions for a planar $4$-web to be of rank $1\;$or$%
\;2\;$or$\;3$''. Finally, we find invariant characterization of planar $5$%
-webs of maximum rank and prove than in general such webs are not
linearizable.
\end{abstract}

\section{Introduction}

Bol in \cite{bo 32} (see also \cite{BB 38} and \cite{B 55}) proved that the
rank of a planar $d$-web does not exceed $(d-1)\left( d-2\right) /2.\;$Chern
in \cite{C82} posed the problem: ``Determine all $d$-webs of curves in the
plane having maximum rank $(d-1)(d-2)/2$, $d$ $\geqq 5.$''

In the current paper, we find an invariant characterization of planar $d$%
-webs of maximum rank and provide a detailed description for the cases $%
d=4,5 $. This\ is the first step for solution of Chern's problem formulated
above.

For $4$-webs, it is well known that the geometry of a $4$-web$\;$is
determined by the curvature $K\;$of one of its $3$-subwebs, the basic
invariant $a$ and their (covariant) derivatives.

We present the characterization of $4$-webs$\;$of maximum rank in two forms,
an invariant analytic form: \emph{A planar }$4$\emph{-web is of maximum rank
three if and only if the curvature }$K\;$\emph{of one of its }$3$\emph{%
-subwebs} \emph{and the covariant derivatives }$K_{3}\;$\emph{and }$K_{4}\;$%
\emph{of }$K\;$\emph{are expressed in terms of the }$4$\emph{-web basic
invariant }$a\;$\emph{and the covariant derivatives of }$a\;$\emph{up to the
third order as indicated in formulas\ of Theorem }\ref{4webmaxrank}, and in
a pure geometric form: \emph{A planar }$4$\emph{-web is of maximum rank
three if and only if it}$\;$\emph{is linearizable and its curvature
vanishes\ }$($\emph{Theorem }\ref{4webmaxrank3})\emph{.\ }Note that the
curvature of a $4$-web\ is a weighted sum of curvatures of its four $3$%
-subwebs.

As far as we know, these characterizations are the first intrinsic
descriptions of $4$-webs of maximum rank expressing maximum rank property in
terms of the web invariants$.$

Note that Dou (see \cite{dou 53} and \cite{dou 55})$\;$studied the rank
problems for planar $4$-webs. The conditions which he found were neither
invariant nor effective. This was a reason that Blaschke (who was familiar
with Dou's results) in his book \cite{B 55} (see $\S 48$, problem $A_{2}$)
listed as open the following problem: ``Find invariant conditions for a
planar $4$-web to be of rank $1\;$or$\;2\;$or$\;3.$''

Our characterizations of planar $4$-webs of maximum rank indicated above
along with characterizations of planar $4$-webs of rank two and one give a
complete solution of the Blaschke problem. The conditions we found are both
invariant and effective, and we applied them to several examples.

Pantazi \cite{Pa 38} found some necessary and sufficient conditions for a
planar web to be of maximum rank. The paper \cite{Pa 38} was followed by the
papers \cite{Pa 40}\ and \cite{Mi 41}. Recently H\'{e}naut \cite{H 04}\ (who
apparently\ was not familiar with the paper \cite{Pa 38}) associated a
connection with the space of abelian equations admitted by a planar web and
proved that this connection is integrable (i.e., its curvature form
vanishes) if and only if$\;$the web$\;$is of maximum rank (see also \cite{Ri}
and \cite{Ri 05}). Pirio in \cite{P 04b} presented a more detailed
exposition of results of Pantazi in \cite{Pa 38} and\ \cite{Pa 40}\ and Mih%
\u{a}ileanu in \cite{Mi 41}. \ Both characterizations (of Pantazi in Pirio's
interpretation and H\'{e}naut) are not given in terms of the web invariants.

Note also that although a geometric description of planar $4$-webs\ of
maximum rank was known (they are algebraizable, i.e., they are equivalent to
$4$-webs formed by the tangents to an algebraic curve of degree four; see %
\cite{BB 38},\ $\S 27$), their invariant characterization was not known.

Theorem \ref{4webmaxrank3} leads to some interesting results in web
geometry. In particular, Theorem \ref{4webmaxrank3}\ implies immediately the
Poincar\'{e}'s theorem (see Corollary \ref{Poincare}). The classical (Poincar%
\'{e}'s) theorem states: \emph{A planar}$\;4$\emph{-web of maximum rank
three is linearizable. }In our exposition, this theorem becomes obvious
because the linearizability conditions are a part of maximum rank conditions
(Theorem \ref{4webmaxrank3}). The Poincar\'{e} theorem was noted in the
books \cite{BB 38} ($\S 27$, p. 239) and \cite{B 55} ($\S 44$) . Note that
this theorem is called the Poincar\'{e} theorem because it is related to
Poincar\'{e}'s mapping (see \cite{Po01}) which is widely used in rank
problems for webs (see, for example, \cite{CG78a} and \cite{CG78b}). It is
worth to note that this theorem can be considered as a reformulation of
Sophus Lie's result on surfaces of double translation (see \cite{Lie82}) in
the web terms. In fact, in web terms, \ Lie's result in \cite{Lie82} means
that a planar $4$-web of maximum rank three is algebraizable (i.e., it is
formed by the tangents to a plane algebraic curve of degree four). This
implies that any $4$-web of maximum rank is linearizable (cf., for example, %
\cite{B 55}, $\S 44$).

Remark also that our proof of the Poincar\'{e} theorem uses essentially the
linearizability conditions found recently in \cite{AGL}.

In this paper, we also find invariant descriptions of $4$-webs of rank two
or one (Theorems \ref{rank2} and \ref{rank1}) and prove that in general such
webs are not linearizable (Propositions \ref{nonlin of webs with 2 abeq} and %
\ref{nonlinear webs of rank 1}). Using theorem \ref{4webmaxrank3}, we prove
also that for linearizable $4$-webs the vanishing of its curvature is not
only necessary but also sufficient for being of maximum rank and that
parallelizable $4$-webs as well as Mayrhoffer's webs are of maximum rank
three.\

We also consider concrete examples of $4$-webs (Examples \ref{4pencils}, \ref%
{3pencils+conics}, \ref{rank2nonlin1}, \ref{rank2nonlin2} and \ref%
{rank1nonlin1})$\;$and applying Theorem \ref{4webmaxrank3},\ establish that
two of them (Examples \ref{4pencils} and Examples \ref{3pencils+conics}) are
of maximum rank, two others (Examples \ref{rank2nonlin1}$\;$and \ref%
{rank2nonlin2})$\;$are of rank two and the last one (\ref{rank1nonlin1}) is
of rank one. Because $4$-webs\ of Examples \ref{rank2nonlin1}, \ref%
{rank2nonlin2}$\;$and \ref{rank1nonlin1} are not linearizable, in general, $%
4 $-webs of ranks two and one are not linearizable. We also study rank
problems for planar $4$- and $5$-webs with constant basic invariants.

\section{Basic constructions for planar webs}

\subsection{Planar $\mathbf{\emph{d}}$-webs}

A $d$-web $W_{d},$ $d\geq 3,$ on a domain $\mathbb{D\subset }\mathbb{R}^{2}$
is defined by $d$ one-dimensional foliations in general position (i.e.,
leaves of any pair of foliations are transversal to each other). Such
foliations can be defined by $d$ functions ($1$-st integrals of the
foliations) $\left\langle f_{1},...,f_{d}\right\rangle $ such that any pair
of functions $f_{i},f_{j},$ $i\neq j,$ are independent, or by $d$
differential \thinspace $1$-forms $<\omega _{1},\omega _{2},\omega
_{3},\omega _{4},...,$ $\omega _{d}>$ such that any two of them are linearly
independent.

We fix a co-basis $\left\langle \omega _{1},\omega _{2}\right\rangle $ and a
$3$-subweb $W_{3}=\left\langle \omega _{1},\omega _{2},\omega
_{3}\right\rangle $. The forms $\omega _{1},\omega _{2},$ and $\omega _{3}$
can be normalized in such a way that
\begin{equation*}
\omega _{1}+\omega _{2}+\omega _{3}=0.
\end{equation*}%
One can easily prove that in this case there is a unique differential $1$%
-form $\gamma $ such that the so-called \emph{structure equations}
\begin{equation*}
d\omega _{i}=\omega _{i}\wedge \gamma
\end{equation*}%
hold for all $i=1,2,3\;$(see \cite{AGL})$.$

The form $\gamma $ determines the Chern connection $\Gamma $ in the
cotangent bundle $T^{\ast }M$ with the following covariant differential:
\begin{equation*}
d_{\Gamma }:\omega _{i}\longmapsto -\omega _{i}\otimes \gamma .
\end{equation*}%
The curvature of this connection is equal to
\begin{equation*}
R_{\Gamma }:\omega _{i}\longmapsto -\omega _{i}\otimes d\gamma .
\end{equation*}%
If we write
\begin{equation*}
d\gamma =K\omega _{1}\wedge \omega _{2},
\end{equation*}%
then the function $K$ is called the \emph{curvature function }of the $3$-web
$W_{3}.$

Note that the curvature form $d\gamma $ is an invariant of the $3$-web $%
W_{3} $ while the curvature function $K$ is a relative invariant of the web.

The scale transformation $\left\langle \omega _{1},\omega _{2},\omega
_{3}\right\rangle \longmapsto \left\langle \omega _{1}^{s},\omega
_{2}^{s},\omega _{3}^{s}\right\rangle ,$ where $s$ is a nonvanishing smooth
function and $\omega _{i}^{s}=s^{-1}\omega _{i},$ preserves the $3$-web in
the sense that triples $\left\langle \omega _{1},\omega _{2},\omega
_{3}\right\rangle $ and $\left\langle \omega _{1}^{s},\omega _{2}^{s},\omega
_{3}^{s}\right\rangle $ determine the same web. The structure equations for $%
\left\langle \omega _{1}^{s},\omega _{2}^{s},\omega _{3}^{s}\right\rangle $
have the form
\begin{equation*}
d\omega _{i}^{s}=\omega _{i}^{s}\wedge \gamma ^{s}
\end{equation*}%
with $\gamma ^{s}=\gamma +d\ln \left\vert s\right\vert ,$ and therefore $%
d\gamma =d\gamma ^{s}.$

If one defines the curvature function $K^{s}$ by the equation
\begin{equation*}
d\gamma ^{s}=K^{s}\omega _{1}^{s}\wedge \omega _{2}^{s},
\end{equation*}%
then%
\begin{equation*}
K^{s}=s^{2}K.
\end{equation*}%
We emphasize this by saying that $K$ is a relative invariant of weight two.

Let $\left\langle \partial _{1},\partial _{2}\right\rangle $ be the basis
dual to $\left\langle \omega _{1},\omega _{2}\right\rangle .$ We put $%
\partial _{3}=\partial _{2}-\partial _{1}.$ Then leaves of the $3$-web $%
W_{3} $ are trajectories of the vector fields $\partial _{2},\partial _{1},$%
and $\partial _{3}.$

We denote by $\delta _{i}$ the covariant derivatives in the direction $%
\partial _{i}$ with respect to the Chern connection.

Let%
\begin{equation*}
\gamma =g_{1}\omega _{1}+g_{2}\omega _{2}.
\end{equation*}%
Then
\begin{equation*}
K=\partial _{1}\left( g_{2}\right) -\partial _{2}\left( g_{1}\right) ,
\end{equation*}%
and the action of the covariant derivatives $\delta _{i}$ on functions $u$
of weight $w$ is:
\begin{equation*}
\delta _{i}^{\left( w\right) }\left( u\right) =\partial _{i}\left( u\right)
-wg_{i}u.
\end{equation*}%
In what follows, we shall skip the superscript when the weight of $u$ is
known. Remark that the covariant derivatives satisfy the Leibnitz rule and
\begin{equation*}
\delta _{2}^{\left( w+1\right) }\delta _{1}^{\left( w\right) }-\delta
_{1}^{\left( w+1\right) }\delta _{2}^{\left( w\right) }=wK
\end{equation*}%
(see \cite{GL 06}).

For general $d$-web $W_{d}=\left\langle \omega _{1},\omega _{2},\omega
_{3},....,\omega _{d}\right\rangle ,$ we choose $\omega _{i}$ for $i\geq 4$
in such a way that the normalizations
\begin{equation*}
a_{i}\omega _{1}+\omega _{2}+\omega _{i+2}=0
\end{equation*}%
hold for $i=1,...,d-2,$ and $a_{1}=1.$

Note that $a_{i}\neq 0,1$ for $i\geq 2.$ Moreover, for the fixed $i,$\ the
value $a_{i}\left( x\right) ,$ $x\in \mathbb{D},$ of the function $a_{i}$ is
the cross-ratio of the four straight lines in $T_{x}^{\ast }(\mathbb{D})$ \
generated by the covectors $\omega _{1,x},\omega _{2,x},\omega _{3,x}$, and $%
\omega _{i+2,x},$ and therefore it is an invariant. The functions $a_{i}$
are called the \emph{basic invariants }(cf. \cite{G 04} or \cite{G 88}, pp.
302--303).

\subsection{Web functions}

We choose (local) coordinates $x,y$ in $\mathbb{D}$ in such a way that $%
\omega _{1}\wedge dx=0\;$and $\omega _{2}\wedge dy=0.$ Let $\omega
_{3}\wedge df$ $=0,$ $\omega _{i+3}\wedge dg_{i}$ $=0,$ $i=1,...,d-3,$ for
some functions $f\left( x,y\right) ,g_{i}\left( x,y\right) .$

Using the scale transformation, we assume that $\omega _{3}=df.$ Then $%
\omega _{1}=-f_{x}dx$ and $\omega _{2}=-f_{y}dy.$

The dual basis $\left\{ \partial _{1},\partial _{2}\right\} $ has the form
\begin{equation*}
\partial _{1}=-f_{x}^{-1}\partial _{x},\ \ \ \partial
_{2}=-f_{y}^{-1}\partial _{y}.
\end{equation*}%
The connection form is
\begin{equation*}
\gamma =-H\omega _{3},
\end{equation*}%
where
\begin{equation*}
g_{1}=g_{2}=H=\frac{f_{xy}}{f_{x}~f_{y}}
\end{equation*}%
(see \cite{GL 06}). The curvature function has the following expression:
\begin{equation*}
K=-f_{x}^{-1}f_{y}^{-1}\left( \log \left( f_{x}f_{y}^{-1}\right) \right)
_{xy}.
\end{equation*}%
In terms of the web functions, the basic invariants have the form%
\begin{equation*}
a_{i}=\frac{f_{y}g_{i+1,x}}{f_{x}g_{i+1,y}}
\end{equation*}%
for $i=2,...,d-2.$

\begin{definition}
A planar $d$-web $W_{d}$ is said to be $($locally$)$ \textbf{parallelizable}
if it is $($locally$)\;$equivalent to a $d$-web of parallel straight lines
in a domain of the affine plane $\mathbb{A}^{2}$.
\end{definition}

It is well known (see, for example, \cite{B 55}, $\S 8$) that \emph{a planar
}$3$\emph{-web is locally parallelizable if and only if }$K=0.$

For planar $d$-webs, $d\geq 4,$ the following statement holds (cf. \cite{G
04} or \cite{G 88}, Section 7.2.1 for $d=4$).

\begin{theorem}
\label{4webpar}A planar $d$-web $W_{d}$ $=\left\langle \omega _{1},\omega
_{2},\omega _{3},\omega _{4},...,\omega _{d}\right\rangle $ is locally
parallelizable if and only if its $3$-subweb $W_{3}$ $=\left\langle \omega
_{1},\omega _{2},\omega _{3}\right\rangle $ is locally parallelizable $($%
i.e., $K=0)$, and all basic invariants $a_{i}$ are constants.
\end{theorem}

\begin{proof}
Let $K=0$ and $a_{i}=\limfunc{const}.$ Then $W_{3}$ $=\left\langle \omega
_{1},\omega _{2},\omega _{3}\right\rangle $ is locally parallelizable, and
we can choose local coordinates $x,y$ in such a way that $\omega
_{1}=-dx,\;\omega _{2}=-dy,$ $\omega _{3}=d(x+y).$ Since $a_{i}=\limfunc{%
const},$ then $\omega _{i+2}=d(a_{i}x+y).$

Conversely,\ suppose that $W_{d}$ $=\left\langle \omega _{1},\omega
_{2},\omega _{3},\omega _{4},...,\omega _{d}\right\rangle $ is locally
parallelizable. We choose local coordinates $x,y$ in such a way that leaves
of the foliations are parallel straight lines in these coordinates. Then
\begin{equation*}
\frac{f_{x}}{f_{y}}\ \text{and\ }\frac{g_{i+1,x}}{g_{i+1,y}}
\end{equation*}%
are constants.

Therefore $K=0,$ and $a_{i}=\limfunc{const}$ due to the above formulae for $%
K $ and $a_{i}.$
\end{proof}

\section{Abelian equations}

\subsection{Classical abelian relations}

We begin with an interpretation of the classical Abel addition theorem (see %
\cite{AB}) in terms of planar webs (cf. \cite{B 55}). A straight line on the
affine plane is defined by a pair $\left( r,s\right) $: $rx+sy=1.$ Assume
that $\left( r,s\right) $ satisfies a cubic equation, say, $%
s^{2}-4r^{3}-g_{2}r-g_{3}=0.$ Given $\left( x,y\right) ,$ one gets the cubic
equation for $r$ of the form $r^{3}+ar^{2}+br+c=0$ with
\begin{equation*}
a=-\frac{x^{2}}{4y^{2}},\,b=\frac{1}{4}\left( g_{2}+\frac{2x}{y}\right) ,c=%
\frac{1}{4}\left( g_{3}-\frac{1}{y^{2}}\right) .
\end{equation*}%
Then in the domain, where
\begin{equation*}
x^{4}-24xy^{2}-12g_{2}y^{4}>0,\ y\neq 0,
\end{equation*}%
the cubic equation has three distinct real roots and consequently three
pairwise independent straight lines $\left( r_{1}\left( x,y\right)
,s_{1}(x,y)\right) ,\left( r_{2}\left( x,y\right) ,s_{2}(x,y)\right) $\ and$%
\ \left( r_{3}\left( x,y\right) ,s_{3}(x,y)\right) $ passing through the
point $\left( x,y\right) .$ They generate a $3$-web $W_{3}$ in the domain.

Let $g_{2}^{3}-27g_{3}^{2}\neq 0.$ Then the solutions of the equation $%
s^{2}-4r^{3}-g_{2}r-g_{3}=0$ can be parametrized by the Weierstrass function
$\wp :\ r=\wp \left( t\right) ,\ s=\wp ^{\prime }\left( t\right) .$ As a
result, the roots $\left( r_{1}\left( x,y\right) ,s_{1}(x,y)\right) ,\left(
r_{2}\left( x,y\right) ,s_{2}(x,y)\right) $\ and$\ \left( r_{3}\left(
x,y\right) ,s_{3}(x,y)\right) $ correspond to three solutions $\left(
t_{1}\left( x,y\right) ,t_{2}\left( x,y\right) ,t_{3}\left( x,y\right)
\right) $ of the equation
\begin{equation*}
f\left( t\right) =\wp \left( t\right) x+\wp ^{\prime }\left( t\right) y-1=0.
\end{equation*}%
Computing the integral
\begin{equation*}
\int t\frac{f^{\prime }\left( t\right) }{f\left( t\right) }dt
\end{equation*}%
along the boundary of the period parallelogram, one finds the Abel relation%
\begin{equation*}
t_{1}\left( x,y\right) +t_{2}\left( x,y\right) +t_{3}\left( x,y\right) =%
\limfunc{const}.
\end{equation*}%
By the construction, the functions $t_{1}\left( x,y\right) ,t_{2}\left(
x,y\right) ,\ $and$\ t_{3}\left( x,y\right) $ are constant on the
corresponding leaves of $W_{3}.$

Consider now an arbitrary planar $d$-web defined by $d$ functions $%
W_{d}=\left\langle f_{1},...,f_{d}\right\rangle .$ Then an \emph{abelian
relation} is given by $d$ functions $(F_{1},....,F_{d})$ of one variable
such that
\begin{equation*}
F_{1}\left( f_{1}\right) +\cdots +F_{d}\left( f_{d}\right) =\limfunc{const}.
\end{equation*}%
We say that two abelian relations $\left( F_{1},....,F_{d}\right) $ and $%
\left( G_{1},....,G_{d}\right) $ are \emph{equivalent }if and only if $%
F_{i}=G_{i}+\limfunc{const}_{i}$ for all $i=1,...,d.$

Obviously the set of equivalence classes of abelian relations admits the
vector space structure with respect to addition: $\left(
F_{1},....,F_{d}\right) +\left( G_{1},....,G_{d}\right) =\left(
F_{1}+G_{1},....,F_{d}+G_{d}\right) $ and multiplication by numbers: $\alpha
\left( F_{1},....,F_{d}\right) =\left( \alpha F_{1},....,\alpha F_{d}\right)
.$ The dimension of this vector space is called the \emph{rank} of the web.

In the case when $d$-web is defined by differential $1$-forms $%
W_{d}=\left\langle \omega _{1},....,\omega _{d}\right\rangle ,$ the
differentiation of the abelian relation leads us to the \emph{abelian
equation }%
\begin{equation*}
\lambda _{1}\omega _{1}+\cdots +\lambda _{d}\omega _{d}=0,
\end{equation*}%
for functions $\left( \lambda _{1},...,\lambda _{d}\right) $ under the
condition that all differential $1$-forms $\lambda _{i}\omega _{i}$ are
closed. The abelian equation is a system of the first order linear PDEs for
the functions $\left( \lambda _{1},...,\lambda _{d}\right) ,$ and the rank
of the web is the dimension of the solution space.

The following example of the $3$-web illustrates the above constructions.
Consider the $3$-web $W_{3}$ given by the web function
\begin{equation*}
f=\frac{2xy-x+y}{x+y}.
\end{equation*}%
Then%
\begin{equation*}
\omega _{1}=-f_{x}dx,\ \omega _{2}=-f_{y}dy,\ \omega _{3}=df.
\end{equation*}%
The condition
\begin{equation*}
\lambda _{1}\omega _{1}+\lambda _{2}\omega _{2}+\lambda _{3}\omega _{3}=0
\end{equation*}%
implies
\begin{equation*}
\lambda _{1}=\lambda _{2}=\lambda _{3},
\end{equation*}%
and the condition $d\left( \lambda _{3}\omega _{3}\right) =0$ gives%
\begin{equation*}
\lambda _{1}=\lambda _{2}=\lambda _{3}=\lambda \left( t\right)
\end{equation*}%
for some function $\lambda (t).$

The other two conditions $d\left( \lambda \omega _{1}\right) =d\left(
\lambda \omega _{2}\right) =0$ lead to the differential equation on $\lambda
$:%
\begin{equation*}
2t\lambda (t)+\left( t^{2}-1\right) \lambda ^{\prime }\left( t\right) =0.
\end{equation*}%
Thus
\begin{equation*}
\lambda \left( t\right) =\frac{1}{t^{2}-1},
\end{equation*}%
and the abelian relation
\begin{equation*}
F_{1}\left( x\right) +F_{2}\left( y\right) +F_{3}\left( f\right) =0
\end{equation*}%
corresponds to the following functions%
\begin{equation*}
F_{1}\left( x\right) =\ln \frac{x+1}{x},\;F_{2}\left( y\right) =\ln \frac{y}{%
y-1},\;F_{3}\left( f\right) =\ln \frac{f-1}{f+1}.
\end{equation*}

\subsection{Abelian differential equations}

In this section we formalize the above constructions. Let $%
W_{d}=\left\langle \omega _{1},...,\omega _{d}\right\rangle $ be a planar $d$%
-web in a domain $\mathbb{D}\subset \mathbb{R}^{2},\;$and let $\pi $:$%
\;E\rightarrow \mathbb{D}$ be a subbundle of the trivial bundle $\mathbb{R}%
^{d}\times \mathbb{D}\rightarrow \mathbb{D}$ consisting of points $\left(
x_{1},...,x_{d},a\right) ,\;$where $\left( x_{1},...,x_{d}\right) \in
\mathbb{R}^{d}$, $a\in \mathbb{D}$, such that $\sum_{1}^{d}x_{i}\omega
_{i,a}=0.$

By the \textbf{abelian equation}\emph{\ }associated with the $d$-web $W_{d}$
we mean a system of first order differential equations for sections $\left(
\lambda _{1},...,\lambda _{d}\right) $ (i.e., $\sum_{1}^{d}\lambda
_{i}\omega _{i}=0$) of the bundle $\pi \ $such that (cf. $\cite{Gr 77}$):%
\begin{equation*}
d\left( \lambda _{1}\omega _{1}\right) =\cdots =d\left( \lambda _{d}\omega
_{d}\right) =0.
\end{equation*}

Let us write down the abelian equation in the explicit form.

\emph{In what follows, we shall choose a }$3$\emph{-subweb, say }$%
\left\langle \omega _{1},\omega _{2},\omega _{3}\right\rangle ,$\emph{\ and
normalize the }$d$\emph{-web as it was done earlier}:%
\begin{equation*}
a_{1}\omega _{1}+\omega _{2}+\omega _{3}=0,\ \;a_{2}\omega _{1}+\omega
_{2}+\omega _{4}=0,....,\;a_{d-2}\omega _{1}+\omega _{2}+\omega _{d}=0,
\end{equation*}%
\emph{with }$a_{1}=1\;$\emph{and} $\ d\omega _{3}=0.$

We call such a normalization \emph{standard}.

Then%
\begin{eqnarray*}
d\left( \lambda _{1}\omega _{1}\right) &=&\left( -\partial _{2}\left(
\lambda _{1}\right) +H\lambda _{1}\right) \omega _{1}\wedge \omega _{2}, \\
d\left( \lambda _{2}\omega _{1}\right) &=&\left( \partial _{1}\left( \lambda
_{2}\right) -H\lambda _{2}\right) \omega _{1}\wedge \omega _{2}, \\
d\left( \lambda _{3}\omega _{3}\right) &=&\left( \partial _{2}\left( \lambda
_{3}\right) -\partial _{1}(\lambda _{3})\right) \omega _{1}\wedge \omega
_{2}, \\
d\left( \lambda _{i}\omega _{i}\right) &=&\left( a_{i-2}\partial _{2}\left(
\lambda _{i}\right) -\partial _{1}(\lambda _{i})+\lambda _{i}\left(
H+\partial _{2}\left( a_{i-2}\right) -a_{i-2}H\right) \right) \omega
_{1}\wedge \omega _{2},
\end{eqnarray*}%
for all $i=4,...,d.$

We shall assume that $\lambda _{i}$ are functions of weight $1$ and $a_{i}$
are of weight $0.$ Then the above formulae take the form%
\begin{eqnarray*}
d\left( \lambda _{1}\omega _{1}\right) &=&-\delta _{2}\left( \lambda
_{1}\right) \omega _{1}\wedge \omega _{2},\  \\
d\left( \lambda _{2}\omega _{1}\right) &=&\delta _{1}\left( \lambda
_{2}\right) \ \omega _{1}\wedge \omega _{2}, \\
d\left( \lambda _{3}\omega _{3}\right) &=&\left( \delta _{2}\left( \lambda
_{3}\right) -\delta _{1}(\lambda _{3})\right) \omega _{1}\wedge \omega
_{2},\  \\
d\left( \lambda _{i}\omega _{i}\right) &=&\left( \delta _{2}\left(
a_{i-2}\lambda _{i}\right) -\delta _{1}(\lambda _{i})\right) \omega
_{1}\wedge \omega _{2}.
\end{eqnarray*}

The normalization condition $\sum_{1}^{d}\lambda _{i}\omega _{i}=0$ implies
that
\begin{eqnarray*}
\lambda _{1} &=&\sum_{1}^{d-2}a_{i}u_{i},\ \lambda
_{2}=\sum_{1}^{d-2}u_{i},\  \\
\lambda _{i+2} &=&u_{i},\ i=1,...,d-2.
\end{eqnarray*}%
Therefore the abelian equation \ is equivalent to the following PDEs system
\begin{align*}
& \Delta _{1}\left( u_{1}\right) =\cdots =\Delta _{d-2}\left( u_{d-2}\right)
=0, \\
& \delta _{1}\left( u_{1}\right) +\cdots +\delta _{1}\left( u_{d-2}\right)
=0,
\end{align*}%
where $\Delta _{i}=\delta _{1}-\delta _{2}\circ a_{i}.$

Let $\mathfrak{A}_{1}\subset \mathbf{J}^{1}\left( \pi \right) $ be the
subbundle of the $1$-jet bundle corresponding to the abelian equation, and $%
\mathfrak{A}_{k}\subset \mathbf{J}^{k}\left( \pi \right) $ be the $\left(
k-1\right) $-prolongation of $\mathfrak{A}_{1}.$ Denote by $\pi _{k,k-1}$:$\;%
\mathfrak{A}_{k}\rightarrow \mathfrak{A}_{k-1}$ the restrictions of the
natural projections $\mathbf{J}^{k}\left( \pi \right) \rightarrow \mathbf{J}%
^{k-1}\left( \pi \right) .$

\begin{proposition}
Let $k\leq d-2.$ Then $\mathfrak{A}_{k}$ are vector bundles and the maps $%
\pi _{k,k-1}$:$\;\mathfrak{A}_{k}\leftarrow \mathfrak{A}_{k-1}$ are
projections. Moreover, $\dim \ker \pi _{k,k-1}=d-k-2.$
\end{proposition}

\begin{proof}
Let $u_{i,r_{1}...r_{s}}$ be coordinates in the jet space $\mathbf{J}%
^{k}\left( \pi \right) $ which correspond to the covariant derivatives $%
\delta _{r_{1}}\cdots \delta _{r_{s}}$ (see \cite{GL 06}\ for more details).
In these coordinates, the abelian equation takes the following form:
\begin{eqnarray*}
&&u_{1,1}=a_{1}u_{1,2}+a_{1,2}u_{1}, \\
&&\cdots \cdots \cdots \cdots \cdots \cdots \cdots \\
&&u_{d-2,1}=a_{d-2}u_{d-2,2}+a_{d-2,2}u_{d-2}, \\
&&u_{1,1}+\cdots +u_{d-2,1}=0.
\end{eqnarray*}%
This means that $u_{1},...,u_{d-2}$ are fiberwise coordinates in the bundle $%
\mathbb{D}\overset{\pi }{\longleftarrow }E,$ while $u_{1,2},...,u_{d-3,2}$
are fiberwise coordinates in the bundle $E\overset{\pi _{1,0}}{%
\longleftarrow }\mathfrak{A}_{1}.$ Taking covariant derivatives of the
abelian equation, we observe that $u_{1,22},...,u_{d-4,22}$ are fiberwise
coordinates in the bundle $\mathfrak{A}_{1}\overset{\pi _{2,1}}{%
\longleftarrow }\mathfrak{A}_{2}$ , etc. This process proves the proposition.
\end{proof}

Proposition 3 shows that there is the following tower of vector bundles:%
\begin{equation*}
\mathbb{D}\overset{\pi }{\longleftarrow }E\overset{\pi _{1,0}}{%
\longleftarrow }\mathfrak{A}_{1}\overset{\pi _{2,1}}{\longleftarrow }%
\mathfrak{A}_{2}\overset{\pi _{3,1}}{\longleftarrow }\cdots \overset{\pi
_{d-3,d-4}}{\longleftarrow }\mathfrak{A}_{d-3}\overset{\pi _{d-2,d-3}}{%
\longleftarrow }\mathfrak{A}_{d-2}.
\end{equation*}%
The last projection $\mathfrak{A}_{d-2}\overset{\pi _{d-2,d-3}}{%
\longrightarrow }\mathfrak{A}_{d-3}$ is an isomorphism, and geometrically it
can be viewed as a linear connection in the vector bundle $\pi _{d-3}$:$\;%
\mathfrak{A}_{d-3}\rightarrow \mathbb{D}.$ Remark that \emph{the abelian
equation is formally integrable if and only if this linear connection is flat%
}.

The dimension of this bundle is equal to $(d-2)+\left( d-3\right) +\cdots
+1=(d-2)\left( d-1\right) /2.$ This shows that the solution space \emph{Sol\
}$\left( \mathfrak{A}\right) $ of the abelian equation\emph{\ }$\mathfrak{A}$
is finite-dimensional and $\dim \;$\emph{Sol\ }$\left( \mathfrak{A}\right)
\leq (d-1)\left( d-2\right) /2.$

The dimension $\dim \;$\emph{Sol\ }$\left( \mathfrak{A}\right) $ is called
the \emph{rank of } the corresponding $d$-web $W_{d}$.

As a consequence, we get the following result which was first established by
Bol \cite{bo 32} (see also \cite{BB 38} and \cite{B 55}).

\begin{theorem}
The rank of a planar $d$-web $W_{d}\;$does not exceed $(d-1)\left(
d-2\right) /2.\;$
\end{theorem}

Remark also, that a different approach for description of the bundle $\pi
_{d-3}$:$\;\mathfrak{A}_{d-3}\rightarrow \mathbb{D}$ in the category of
analytical webs was used in \cite{H 04}.

The obstruction for compatibility of the abelian equation is given by the
multi-bracket (see \cite{KL 06} and Section 7.1). The matrix of the abelian
system is

\begin{equation*}
\begin{Vmatrix}
\Delta _{1} & \cdots & 0 \\
\vdots & \ddots & \vdots \\
0 & \cdots & \Delta _{d-2} \\
\delta _{1} & \cdots & \delta _{1}%
\end{Vmatrix}%
.
\end{equation*}

Computing the multi-bracket, we get%
\begin{eqnarray*}
&&\left( -1\right) ^{d}\left\{ \left( \Delta _{1},..,0\right) ,....,\left(
0,...,\Delta _{d-2}\right) ,\left( \delta _{1},...,\delta _{1}\right)
\right\} =\delta _{1}\Delta _{2}\cdots \Delta _{d-2}\left( \Delta
_{1},..,0\right) \\
&&+\Delta _{1}\delta _{1}\Delta _{3}\cdots \Delta _{d-2}\left( 0,\Delta
_{2},..,0\right) +\cdots +\Delta _{1}\cdots \Delta _{i-1}\delta _{1}\Delta
_{i+1}\cdots \Delta _{d-2}\left( 0,...,\Delta _{i},...,0\right) \\
&&+\cdots +\Delta _{1}\cdots \Delta _{d-3}\delta _{1}\left( 0,...,\Delta
_{d-2}\right) -\Delta _{1}\cdots \Delta _{d-2}\left( \delta _{1},...,\delta
_{1}\right) .
\end{eqnarray*}%
Therefore, the compatibility condition for the abelian system is%
\begin{equation*}
\varkappa =\square _{1}u_{1}+\cdots +\square _{d-2}u_{d-2}=0,
\end{equation*}%
where
\begin{equation*}
\square _{i}=\Delta _{1}\cdots \Delta _{d-2}\cdot \delta _{1}-\Delta
_{1}\cdots \Delta _{i-1}\cdot \delta _{1}\cdot \Delta _{i+1}\cdots \Delta
_{d-2}\cdot \Delta _{i}
\end{equation*}%
are linear differential operators of order not exceeding$\;d-2.$

Summarizing, we get the following

\begin{theorem}
A $d$-web is of maximum rank $(d-1)\left( d-2\right) /2$ if and only if $%
\varkappa =0$ on $\mathfrak{A}_{d-2}.$
\end{theorem}

Remark that $\varkappa $ can be viewed as a linear function on the vector
bundle $\mathfrak{A}_{d-2},$ and therefore the above theorem imposes $%
(d-1)\left( d-2\right) /2$ conditions on the $d$-web (or on $d-2$ web
functions) in order the web has the maximum rank. A calculation of these
conditions is pure algebraic, and we shall illustrate this calculation below
for planar $\;3$-,\ $4$- and $5$-webs. All these calculations are based on
expressions for total covariant derivatives given in \cite{GL 06}. Note also
that expressions for $\varkappa $ in the case of general $d$-webs are
extremely cumbersome while for concrete $d$-webs it is not the case.

\section{Rank of a planar $\mathbf{3}$-web}

Let $d=3.\;$Then \ the maximum rank of $W_{3}\;$is $1.$The abelian equation
takes the form%
\begin{equation*}
\begin{array}{l}
\Delta _{1}\left( u_{1}\right) =0, \\
\delta _{1}\left( u_{1}\right) =0.%
\end{array}%
\end{equation*}%
The obstruction $\varkappa $ equals
\begin{equation*}
\varkappa =(\delta _{1}\Delta _{1}-\Delta _{1}\delta _{1})u_{1}=\left(
\delta _{2}\delta _{1}-\delta _{1}\delta _{2}\right) u_{1}=Ku_{1}.
\end{equation*}

\begin{theorem}
A $3$-web $W_{3}\;$is of maximum rank one if and only if it is
parallelizable. The only abelian equation admitted by such a $3$-web is the
equation%
\begin{equation*}
\omega _{1}+\omega _{2}+\omega _{3}=0
\end{equation*}%
for the standard normalization.
\end{theorem}

Note that the above theorem was first proved in \cite{BD 28}.\

\section{Planar $\mathbf{4}$-webs}

\subsection{The obstruction}

In the standard normalization for $4$-webs\ $W_{4},$ we put $a_{2}=a$:%
\begin{equation*}
\begin{array}{l}
\omega _{1}+\omega _{2}+\omega _{3} =0,\ \; \\
a\ \omega _{1}+\omega _{2}+\omega _{4} =0%
\end{array}%
\end{equation*}
and reserve the subscripts for the covariant derivatives of $a.$ Thus $%
a_{2}=\delta _{2}\left( a\right) $ and so on.

In what follows, we use the following form of the abelian relation:%
\begin{equation*}
\left( u+av\right) \omega _{1}+\left( u+v\right) \omega _{2}+u\omega
_{3}+v\omega _{4}=0,
\end{equation*}%
where $\lambda _{1}=u+av,\lambda _{2}=u+v,\lambda _{3}=u,\lambda _{4}=v,\ $%
and all summands are closed $1$-forms under condition that $u$ and $v$
satisfy the abelian equation%
\begin{equation*}
\delta _{1}\left( u\right) -\delta _{2}\left( u\right) =0,~\delta _{1}\left(
v\right) -\delta _{2}\left( av\right) =0,~\delta _{1}\left( u\right) +\delta
_{1}\left( v\right) =0.\
\end{equation*}

For $4$-webs,\ the tower of prolongations of the abelian equation is
\begin{equation*}
\mathbb{D}\overset{\pi }{\longleftarrow }E\overset{\pi _{1,0}}{%
\longleftarrow }\mathfrak{A}_{1}\overset{\pi _{2,1}}{\longleftarrow }%
\mathfrak{A}_{2},
\end{equation*}%
where $\pi _{2,1}:\mathfrak{A}_{2}\rightarrow \mathfrak{A}_{1}$ defines a
linear connection on the $3$-dimensional vector bundle $\pi _{1}$:$\;%
\mathfrak{A}_{1}\rightarrow \mathbb{D}.$

We shall use canonical fiberwise coordinates $u,v$ and $u_{1},....,$ etc. in
the jet bundles instead of $u_{1},u_{2},u_{1,1},....$

In these coordinates, the abelian equation takes the form%
\begin{equation*}
u_{1}-u_{2}=0,v_{2}-av_{2}-a_{2}v=0,u_{1}+v_{1}=0,
\end{equation*}%
and the obstruction
\begin{equation*}
\varkappa =(\Delta _{1}\Delta _{2}\delta _{1}-\delta _{1}\Delta _{1}\Delta
_{2})u+(\Delta _{1}\Delta _{2}\delta _{1}-\Delta _{1}\delta _{1}\Delta _{2})v
\end{equation*}%
equals $\varkappa =c_{0}v_{2}+c_{1}v+c_{2}u$.

The straightforward computation gives the following result.

\begin{theorem}
\label{restriction of kappa}In the canonical coordinates, the restriction $%
\varkappa $ on $\mathfrak{A}_{2}$ has the form
\begin{equation}
\varkappa =c_{0}v_{2}+c_{1}v+c_{2}u,  \label{kappa}
\end{equation}%
where%
\begin{eqnarray*}
c_{0} &=&K+\frac{a_{11}-aa_{22}-2\left( 1-a\right) a_{12}}{4a(1-a)}+\frac{%
\left( -1+2a\right) a_{1}^{2}-a^{2}a_{2}^{2}+2\left( 1-a\right)
^{2}a_{1}a_{2}}{4\left( 1-a\right) ^{2}a^{2}}, \\
c_{1} &=&\frac{K_{2}-K_{1}}{4(1-a)}+\frac{\left( a-4\right) a_{1}+\left(
11-20a+12a^{2}\right) a_{2}}{12\left( 1-a\right) ^{2}a}K+\frac{%
a_{112}-a_{122}}{4a(1-a)} \\
&&+\frac{a_{1}-aa_{2}}{4a^{2}(1-a)}a_{22}+\frac{\left( 2a-1\right) \left(
a_{1}-aa_{2}\right) }{4\left( 1-a\right) ^{2}a^{2}}a_{12}-\frac{%
a_{2}^{2}\left( \left( 1-2a\right) a_{1}+aa_{2}\right) }{4\left( 1-a\right)
^{2}a^{2}}, \\
c_{2} &=&\frac{aK_{2}-K_{1}}{4a(1-a)}+\frac{\left( 1-2a\right) a_{1}-\left(
a-2\right) aa_{2}}{4\left( 1-a\right) ^{2}a^{2}}K.
\end{eqnarray*}
\end{theorem}

The coefficient $c_{0}$ in the expression of $\varkappa $ has an intrinsic
geometric meaning. Namely, let us define a curvature of$\ $a $4$-web as the
arithmetic mean of the curvatures of its $3$-subwebs. \ More precisely,
consider the $3$-subwebs $[1,2,3],\;[1,2,4],\;[1,3,4]\;$and $[2,3,4]\;$of a $%
4$-web $W_{4}$ with the following normalizations given by the $1$-forms and
the basic invariant:

\begin{itemize}
\item[{$[1,2,3]$:}] $\omega _{1},\omega _{2},\omega _{3},$ $\omega _{4};$

\item[{$[1,2,4]$:}] $\rho _{1}=a\omega _{1},\rho _{2}=\omega _{2},\rho
_{3}=\omega _{4},\rho _{4}=\omega _{3};$

\item[{$[1,3,4]$:}] $\sigma _{1}=(a-1)\omega _{1},\sigma _{2}=-\omega
_{3},\sigma _{3}=\omega _{4},\sigma _{4}=-\omega _{2};$

\item[{$[2,3,4]$:}] $\tau _{1}=(a-1)\omega _{2},\tau _{2}=a\omega _{3},\tau
_{3}=-\omega _{4},\tau _{4}=a\omega _{1}.$
\end{itemize}

Let $K[l,m,n]$ be the curvature function of the $3$-subweb $[l,m,n]$. Define
a \emph{curvature }$2$\emph{-form }$L\omega _{1}\wedge \omega _{2}$\emph{\
of the }$4$\emph{-web} as follows%
\begin{eqnarray*}
4L\omega _{1}\wedge \omega _{2} &=&K[1,2,3]\omega _{1}\wedge \omega
_{2}+K[1,2,4]\rho _{1}\wedge \rho _{2} \\
&&+K[1,3,4]\sigma _{1}\wedge \sigma _{2}+K[2,3,4]\tau _{1}\wedge \tau _{2}.
\end{eqnarray*}%
Then (see \cite{G 04} or Ch. 7 of \cite{G 88} for details)%
\begin{eqnarray*}
K[1,2,3] &=&K, \\
K[1,2,4] &=&\frac{1}{a}\left( K-\frac{a_{12}}{a}+\frac{a_{1}a_{2}}{a^{2}}%
\right) , \\
K[1,3,4] &=&\frac{1}{a-1}\left[ K+\frac{a_{2}(a_{1}-a_{2})}{(1-a)^{2}}+\frac{%
a_{12}-a_{22}}{1-a}\right] , \\
K[2,3,4] &=&\frac{1}{a(a-1)}\left[ K+\frac{(2a-1)a_{1}(a_{1}-a_{2})}{%
a^{2}(1-a)^{2}}+\frac{a_{11}-a_{12}}{a(1-a)}\right] .
\end{eqnarray*}

Computing the \emph{curvature function} $L$ from the above formulae, we
obtain the following geometric interpretations of the coefficient $c_{0}.$

\begin{theorem}
The coefficient $c_{0}$ equals the curvature function of the $4$-web:%
\begin{equation*}
c_{0}=L.
\end{equation*}
\end{theorem}

\subsection{$\mathbf{4}$-webs of maximum rank}

A planar $4$-web has the maximum rank three if and only if the obstruction $%
\varkappa $ identically equals zero, i.e., if and only if $%
c_{0}=c_{1}=c_{2}=0.$ This leads us to the following result.

\begin{theorem}
\label{4webmaxrank}A planar $4$-web $W_{4}$ is of maximum rank three if and
only if its curvature $K\;$and the covariant derivatives $K_{3}\;$and $%
K_{4}\;$of $K,$ where $\partial _{3}=\partial _{2}-\partial _{1}$ and $%
\partial _{4}=a\partial _{2}-\partial _{1},$ are expressed in terms of the $%
4 $-web basic invariant $a\;$and its covariant derivatives up to the third
order as follows:%
\begin{eqnarray*}
K &=&\frac{-a_{11}+aa_{22}+2\left( 1-a\right) a_{12}}{4a(1-a)}+\frac{\left(
1-2a\right) a_{1}^{2}+a^{2}a_{2}^{2}-2\left( 1-a\right) ^{2}a_{1}a_{2}}{%
4\left( 1-a\right) ^{2}a^{2}}, \\
K_{3} &=&\frac{\left( 4-a\right) a_{1}-\left( 11-20a+12a^{2}\right) a_{2}}{%
3\left( 1-a\right) a}K+\frac{a_{122}-a_{112}}{a}+\frac{a_{4}a_{22}}{a^{2}} \\
&&+\frac{\left( 2a-1\right) a_{4}a_{12}}{\left( 1-a\right) a^{2}}+\frac{%
2a_{2}^{2}\left( 1-a\right) a_{1}+a_{2}^{2}a_{4}}{\left( 1-a\right) a^{2}},
\\
K_{4} &=&\frac{aa_{4}-\left( 1-a\right) a_{1}-2aa_{3}}{\left( 1-a\right) a}K.
\end{eqnarray*}
\end{theorem}

Taking the covariant derivatives $\delta _{3}$ and $\delta _{4}$ of the
first equation in the above theorem,\ we find the values of $K_{3}\;$and$%
\;K_{4}.\;$Comparing the obtained values with their values in the theorem,
we arrive at two relations\ (see them below in Proposition \ref{4webmaxrank2}%
) between the $4$-web basic invariant $a\;$and its covariant derivatives up
to the third order.

Conversely, these relations along with the values of $K_{3}\;$and $K_{4}\;$%
obtained by differentiation of $K$ allow us to reconstruct the second and
the third equations of the above theorem.

This proves the following result.

\begin{proposition}
\label{4webmaxrank2}A planar $4$-web is of maximum rank three if and only if
its curvature $K\;$has the form indicated in the first equation of Theorem %
\ref{4webmaxrank}, and the $4$-web basic invariant $a\;$and its covariant
derivatives up to the third order satisfy the following two relations:%
\begin{equation*}
\begin{array}{l}
6(a-1)^{2}a^{2}[-a_{111}+2(a+1)a_{112}-3aa_{122}] \\
+a(a-1)[a(5(7a-5)a_{1}-3(4a^{2}+5a-4)a_{2})a_{11} \\
-2((13a^{2}+18a-19)a_{1}+3a(3-5a)a_{2})a_{12} \\
+a((19a-17)a_{1}+15aa_{2})a_{22} \\
+(-34a^{2}+49a-19)a_{1}^{3}+(26a^{3}+40a^{2}-89a+38)a_{1}^{2}a_{2} \\
+a(-31a^{2}+53a-18)a_{1}a_{2}^{2}-15a^{3}a_{2}^{3} =0%
\end{array}%
\end{equation*}
and
\begin{equation*}
\begin{array}{l}
6(a-1)^{2}a^{2}[3a_{112}-2(a+1)a_{122}+aa_{222}] \\
+a(a-1)[(-15a_{1}+(17-19a)a_{2})a_{11} \\
+2(3(7-9a)a_{1}+(5a^{2}+18a-11)a_{2})a_{12} \\
+(3(4a^{2}+5a-4)a_{1}+a(1-11a)a_{2})a_{22} \\
+15(2a-1)a_{1}^{3}+(56a^{2}-101a+41)a_{1}^{2}a_{2} \\
+(-10a^{3}-41a^{2}+58a-22)a_{1}a_{2}^{2}+a^{2}(5a-1)a_{2}^{3} =0%
\end{array}%
\end{equation*}
\end{proposition}

This proposition allows us to find a geometric meaning of the last two
equations of Theorem \ref{4webmaxrank}.

\begin{proposition}
If the curvature of a $4$-web vanishes,$\;$then conditions of Proposition %
\ref{4webmaxrank2}$\;$are equivalent to linearizability of the $4$-web.
\end{proposition}

\begin{proof}
It is easy to check that under condition $L=0,$ the $4$-web linearizability
conditions given in \cite{AGL} are equivalent to the conditions in
Proposition \ref{4webmaxrank2}.
\end{proof}

Now we can formulate Theorem \ref{4webmaxrank} in pure geometric terms.

\begin{theorem}
\label{4webmaxrank3}A $4$-web$\;$is of maximum rank three if and only if it
is linearizable and its curvature vanishes.
\end{theorem}

\begin{remark}
As far as we know, the above characterizations of $4$-webs of maximum rank
are the first invariant intrinsic descriptions of such webs in terms of the
web invariants.$\;$Moreover, conditions for a $4$-web$\;$to\ be of maximum
rank include the web linearizability conditions.
\end{remark}

Thus, we have three different (but equivalent) invariant analytic conditions
which are necessary and sufficient for a $4$-web$\;$to be of maximum rank
three:

\begin{description}
\item[(i)] The conditions of Theorem \ref{4webmaxrank};

\item[(ii)] Vanishing of the curvature of the $4$-web$\;$and the conditions
of Proposition \ref{4webmaxrank2} ; and

\item[(iii)] Vanishing of the curvature of the $4$-web$\;$and the $4$-web
linearizability conditions from \cite{AGL}.
\end{description}

Each of these three conditions is effective and can be used as a test for
determination whether some concrete $4$-web is of maximum rank (see examples
at the end of this section).

Theorem \ref{4webmaxrank3} leads to some interesting results in web geometry.

For \emph{linearizable} $4$-webs, the condition of vanishing of the curvature%
$\;$is necessary and sufficient for a $4$-web to be of maximum rank.

\begin{corollary}
\label{lin-le max rank 4-web}A linearizable planar$\;4$-web$\;$is \textit{of
maximum rank three if and only if the curvature }vanishes.
\end{corollary}

\textbf{Remark}. The proof of Theorem \ref{4webmaxrank3} (and Corollary \ref%
{lin-le max rank 4-web}) is \ heavily based on the $4$-web linearizability
conditions in \cite{AGL}. For a \emph{linear }$4$-web, the result of
Corollary \ref{lin-le max rank 4-web} was announced (not proved) in \cite{Pa
40} (see also \cite{P 04b}, Section 5.1.3). Our result is more general than
the result for \emph{linear }$4$-webs in \cite{Pa 40}. $\;$

The next corollary gives the direct web-theoretical proof of the Poincar\'{e}
theorem.

\begin{corollary}
\label{Poincare}$(\mathbf{Theorem\;of\;Poincar\acute{e}})\;$A planar$\;4$-web%
$\;$\textit{of maximum rank three} is linearizable.
\end{corollary}

\begin{proof}
The result follows directly from Theorem \ref{4webmaxrank3}, because the
linearizability conditions are a part of conditions of Theorem \ref%
{4webmaxrank3}.
\end{proof}

\begin{corollary}
If a planar $4$-web\ with a constant basic invariant $a\;$has maximum rank
three, then it is parallelizable.
\end{corollary}

\begin{proof}
In fact, if $a=\limfunc{const}.,\;$then $a_{i}=a_{ij}=a_{ijk}=0.\;$If the$%
\;4 $-web is of maximum rank, then by Theorem \ref{4webmaxrank3}, $L=0.$
Substituting $a_{i}=a_{ij}=a_{ijk}=0\;$into $L=0$, we get $K=0.\;$Therefore,
the web is parallelizable.
\end{proof}

\begin{corollary}
Parallelizable planar $4$-webs\ have maximum rank three.
\end{corollary}

\begin{proof}
By Proposition \ref{4webpar}, a $4$-web is parallelizable if and only if \
the following conditions are satisfied:%
\begin{equation*}
K=0,\;\;\;a=\limfunc{const}.
\end{equation*}%
It follows that $L=0$. Because a parallelizable $4$-web\ is linearizable, by
Theorem \ref{4webmaxrank3}, such a web is of maximum rank three.
\end{proof}

\begin{definition}
\label{def of MW}$4$-webs all $3$-subwebs of which are parallelizable $($%
hexagonal$)$ are called \textbf{Mayrhofer }$\mathbf{4}$\textbf{-webs.}
\end{definition}

They were introduced by Mayrhofer (see \cite{M 28}). The following corollary
gives new property of Mayrhofer's $4$-webs.$\;$

\begin{corollary}
\label{MW of rank 3}The Mayrhofer $4$-webs are of maximum rank three.
\end{corollary}

\begin{proof}
First note that by Definition \ref{def of MW}, we have $L=0.\;$Second, it is
well known (see \cite{BB 38}, $\S 10$; see also \cite{G 04}) that the
Mayrhofer $4$-webs are linearizable. Thus, by Theorem \ref{4webmaxrank3},
the Mayrhofer $4$-webs are of, by maximum rank three.
\end{proof}

Note that the result of Corollary \ref{MW of rank 3} was also proved in the
recent paper \cite{Ri}.

In the same way which we used to define the curvature of a $4$-web, by
taking alternative sums, we can find three additional second-order
invariants which are expressed only in terms of the basic invariant $a\;$and
its covariant derivatives of the first and second order:%
\begin{eqnarray*}
M &=&K[1,2,3]-aK[1,2,4]-(a-1)K[1,3,4]+a(a-1)K[2,3,4], \\
P &=&K[1,2,3]+aK[1,2,4]-(a-1)K[1,3,4]-a(a-1)K[2,3,4], \\
Q &=&K[1,2,3]-aK[1,2,4]+(a-1)K[1,3,4]-a(a-1)K[2,3,4].
\end{eqnarray*}%
Then%
\begin{eqnarray*}
M &=&\frac{-a_{11}-2aa_{12}-aa_{22}}{a(a-1)}+\frac{%
(2a-1)a_{1}^{2}-2a^{2}a_{1}a_{2}+a^{2}a_{2}^{2}}{a^{2}(a-1)^{2}}, \\
P &=&\frac{(a_{11}-aa_{22})}{a(a-1)}+\frac{(1-2a)a_{1}^{2}+a^{2}a_{2}^{2}}{%
a^{2}(a-1)^{2}}, \\
Q &=&\frac{a_{11}-2a_{12}+aa_{22}}{a(a-1)}+\frac{%
(1-2a)a_{1}^{2}+2(2a-1)a_{1}a_{2}-a^{2}a_{2}^{2}}{a^{2}(a-1)^{2}}.
\end{eqnarray*}%
$\;$

Using these invariants, we can now establish a new invariant
characterization of \ Mayrhofer's $4$-webs.

\begin{proposition}
A $4$-web is Mayrhofer's web if and only if the invariants $M,P,\newline
Q\;$and $L\; $vanish.
\end{proposition}

\begin{proof}
Consider the system $L=M=P=Q=0\;$as a linear homogeneous system\ with
respect to $K[1,2,3],\;K[1,2,4],\;K[1,3,4]\;$and $K[2,3,4].\;$The
determinant of this system is equal to $-16a^{2}(a-1)^{2}.\;$Because $a\neq
0,1,\;$the determinant is different from $0$. Thus the system has only the
trivial solution $K[1,2,3]=K[1,2,4]=K[1,3,4]=K[2,3,4]=0.\;$

Therefore, by Definition \ref{def of MW}, the $4$-web is a Mayrhofer $4$-web.%
$\;$The converse statement is obvious.
\end{proof}

\subsubsection{Examples}

Remind that we use the following form of the abelian relation:%
\begin{equation*}
\left( u+av\right) \omega _{1}+\left( u+v\right) \omega _{2}+u\omega
_{3}+v\omega _{4}=0,
\end{equation*}%
where all summands are closed $1$-forms under condition that $u$ and $v$
satisfy the abelian equation%
\begin{equation*}
\delta _{1}\left( u\right) -\delta _{2}\left( u\right) =0,~\delta _{1}\left(
v\right) -\delta _{2}\left( av\right) =0,~\delta _{1}\left( u\right) +\delta
_{1}\left( v\right) =0.\
\end{equation*}

The following two cases are important in applications:

\begin{itemize}
\item[$v=0$:] This will be the case if and only if $K=0,$ and the abelian
relation has the form%
\begin{equation*}
u\omega _{1}+u\omega _{2}+u\omega _{3}=0.
\end{equation*}

\item[$u=0$:] In this case the abelian equation gives%
\begin{equation*}
\delta _{1}\left( v\right) =0,\delta _{2}\left( v\right) =\frac{a_{2}v}{a},
\end{equation*}%
and the compatibility conditions%
\begin{eqnarray*}
\delta _{2}\delta _{1}v-\delta _{1}\delta _{2}v &=&Kv, \\
\delta _{2}\delta _{1}v-\delta _{1}\delta _{2}v &=&\delta _{1}\left( \frac{%
a_{2}}{a}\right) v
\end{eqnarray*}%
imply
\begin{equation*}
K=\delta _{1}\left( \frac{a_{2}}{a}\right) =\frac{aa_{12}-a_{1}a_{2}}{a^{2}}.
\end{equation*}%
The abelian relation becomes%
\begin{equation*}
av\omega _{1}+v\omega _{2}+v\omega _{4}=0.
\end{equation*}
\end{itemize}

There are three cases when both these conditions hold, and therefore
\begin{equation*}
u\omega _{1}+u\omega _{2}+u\omega _{3}=0,
\end{equation*}%
and
\begin{equation*}
av\omega _{1}+v\omega _{2}+v\omega _{4}=0
\end{equation*}%
are abelian relations:

\begin{enumerate}
\item Parallelizable $4$-webs;

\item Mayrhofer $4$-webs; and

\item $4$-webs for which $K[1,2,3]=K[1,2,4]=0.$
\end{enumerate}

\begin{example}
\label{4pencils}We consider the planar $4$-web formed by the coordinate
lines $y=\limfunc{const}.,\;x=\limfunc{const}.,\;$and by the level sets of
the functions%
\begin{equation*}
f(x,y)=\frac{x}{y}\;\text{and\ }g(x,y)=\frac{1-y}{1-x}.
\end{equation*}%
\
\end{example}

Note that this is the $4$-subweb of the famous Bol $5$-web which has the
maximum rank six but not linearizable (see Example $7\;$in Section $5.2\;$of %
\cite{AGL}). Note also that the third and the fourth foliations of this $4$%
-web are the pencils of straight lines with the centers at points $(0,0)\;$%
and\ $(1,1).\;$This $4$-web is linear (and therefore linearizable).

First, note that because the $3$-subweb $[1,2,3]$ of this $4$-web is
parallelizable, the web admits the abelian relation$\mathrm{\;}$
\begin{equation*}
u\omega _{1}+u\omega _{2}+u\omega _{3}=0.
\end{equation*}

The direct calculation shows that for this $4$-web the conditions of Theorem %
\ref{4webmaxrank} are satisfied. \ Moreover, the straightforward
computations show that this $4$-web is a Mayrhofer $4$-web, and the latter
is of maximum rank three by Corollary \ref{MW of rank 3}. The corresponding
abelian relations are%
\begin{equation*}
\renewcommand{\arraystretch}{1.5}%
\begin{array}{l}
\ln f_{1}-\ln f_{2}-\ln f_{3}=0, \\
\ln (1-f_{1})-\ln (1-f_{2})+\ln f_{4}=0, \\
\ln \frac{1-f_{1}}{f_{1}}-\ln \frac{1-f_{3}}{f_{3}}-\ln \left(
1-f_{4}\right) =0,%
\end{array}%
\renewcommand{\arraystretch}{1}
\end{equation*}%
where%
\begin{equation*}
f_{1}=x,\;f_{2}=y,\;f_{3}=\frac{x}{y},\text{\ }f_{4}=\frac{1-y}{1-x}.
\end{equation*}

\begin{example}
\label{3pencils+conics}We consider the planar $4$-web formed by the
coordinate lines $y=\limfunc{const}.,\;x=\limfunc{const}.,\;$and by the
level sets of the functions%
\begin{equation*}
f(x,y)=\frac{x}{y}\;\text{and\ }g(x,y)=\frac{x-xy}{y-xy}
\end{equation*}%
$($see Example $8$ in Section $5.2\;$of$\;\cite{AGL})$.\
\end{example}

Note that this is another $4$-subweb of the famous Bol $5$-web. Note also
that the third and the fourth foliations of this $4$-web are the pencil of
straight lines with the center at the point $(0,0)\;$and the foliation of
conics$.\;$It was proved in \cite{AGL} that this $4$-web is linearizable.

By the same reason as in Example \ref{4pencils}, we have again $K=0\;$and
conditions of Theorem \ref{4webmaxrank} are satisfied. Thus, the planar $4$%
-web in question is of maximum rank three with the following abelian
relations:%
\begin{equation*}
\renewcommand{\arraystretch}{1.5}%
\begin{array}{l}
\ln f_{1}-\ln f_{2}-\ln f_{3}=0, \\
\ln \left( \frac{1}{f_{1}}-1\right) -\ln \left( \frac{1}{f_{2}}-1\right)
+\ln f_{4}=0, \\
\ln (1-f_{1})-\ln \left( 1-f_{3}\right) +\ln (1-f_{4})=0,%
\end{array}%
\renewcommand{\arraystretch}{1}
\end{equation*}%
where
\begin{equation*}
f_{1}=x,\;f_{2}=y,\;f_{3}=\frac{x}{y},\text{\ }f_{4}=\frac{x(1-y)}{y(1-x)}.
\end{equation*}

\begin{example}
We consider the planar $4$-web formed by the coordinate lines $y=\limfunc{%
const}.,\;x=\limfunc{const}.,\;$and by the level sets of the functions%
\begin{equation*}
f(x,y)=x+y\;\text{and\ }g(x,y)=x^{2}+y^{2}.
\end{equation*}
\end{example}

By the same reason as in Example \ref{4pencils}, we have again $K=0,\;$and
therefore the web admits the abelian relation$\mathrm{\;}$
\begin{equation*}
u\omega _{1}+u\omega _{2}+u\omega _{3}=0.
\end{equation*}

One can check that the $4$-web linearizability conditions from \cite{AGL}
are not satisfied. Therefore, this $4$-web is not linearizable. By Theorem %
\ref{Poincare},$\;$this $4$-web is not of maximum rank three. Thus, the rank
of the $4$-web in question can be $1$ or $2$. $\;$

\subsection{$\mathbf{4}$-webs of rank two}

As we have seen earlier, a $4$-web admits an abelian equation (has a
positive rank) if and only if the equation
\begin{equation}
c_{0}v_{2}+c_{1}v+c_{2}u=0  \label{eq4}
\end{equation}%
has a nonzero solution.

Suppose that the coefficient $c_{0}\;$in equation (\ref{eq4}) equals $%
0,\;c_{0}=0.\ $Then if two other coefficients $c_{1}\;$and $c_{2}\;$of (\ref%
{eq4}) are also $0,\;$then as we know (see Theorem \ref{4webmaxrank}), a $4$%
-web\ is of maximum rank three.\ If $c_{0}=0\;$but one of the coefficients $%
c_{1}\;$or $c_{2}\;$of (\ref{eq4}) is not $0,\;$then $c_{1}v+c_{2}u=0$ and,
say $u,$ satisfies a $1$-st order PDEs system of two equations. Therefore,
the $4$-web$\;$admits not more than one abelian equation (i.e., it$\;$is of
rank one or zero).

In what follows , we assume that the coefficient $c_{0}\;$in (\ref{eq4}) is
different from $0$:$\;c_{0}\neq 0.\;$Then a $4$-web$\;$cannot be of rank
more that two.

In this section we shall consider the case when \emph{a }$4$\emph{-web is of
rank two}.

\begin{theorem}
\label{rank2}A planar $4$-web is of rank two if and only if $c_{0}\neq 0,$
and
\begin{equation}
G_{ij}=0,i,j=1,2,  \label{cond for rank 2}
\end{equation}%
where%
\begin{eqnarray*}
G_{11} &=&ac_{0}(c_{2,2}-c_{2,1})+ac_{2}(c_{0,1}-c_{0,2})-a\left( 1-a\right)
c_{1}c_{2} \\
&&+\left( 2a_{2}-a_{1}-aa_{2}\right) c_{0}c_{2}-Kc_{0}^{2}, \\
G_{12} &=&ac_{0}(c_{1,2}-c_{1,1})+ac_{1}(c_{0,1}-c_{0,2})-a\left( 1-a\right)
c_{1}^{2} \\
&&+\left( 2a_{2}-a_{1}-2aa_{2}\right) c_{0}c_{1}+\left(
a_{2}^{2}+a_{12}-a_{22}\right) c_{0}^{2}, \\
G_{21}
&=&c_{0}(c_{2,1}-ac_{2,2})+c_{2}(ac_{0,2}-c_{0,1})-2a_{2}c_{0}c_{2}+a\left(
1-a\right) c_{2}^{2}, \\
G_{22} &=&c_{0}(c_{1,1}-ac_{1,2})+c_{1}(ac_{0,2}-c_{0,1})+a\left( 1-a\right)
c_{1}c_{2}-a_{2}c_{0}c_{1} \\
&&-a_{2}(1-a)c_{0}c_{2}+\left( a_{22}-K\right) c_{0}^{2}.
\end{eqnarray*}
\end{theorem}

\begin{proof}
Adding the compatibility condition (\ref{eq4}) to the abelian equations and
solving the resulting system with respect to $u_{1},u_{2},v_{1},\ $and$\
v_{2},\;$we get the Frobenius type PDEs system:%
\begin{eqnarray*}
u_{1} &=&-a_{2}v+\frac{a}{c_{0}}(c_{2}u+c_{1}v), \\
u_{2} &=&-a_{2}v+\frac{a}{c_{0}}(c_{2}u+c_{1}v), \\
v_{1} &=&a_{2}v-\frac{a}{c_{0}}(c_{2}u+c_{1}v), \\
v_{2} &=&-\frac{a}{c_{0}}(c_{2}u+c_{1}v).
\end{eqnarray*}%
We get the integrability conditions for this system from the commutation
relation$\ \delta _{2}\delta _{1}-\delta _{1}\delta _{2}=K.$ Computing the
commutators and substituting $u_{1},u_{2},v_{1}\;$and $v_{2},$ due to the
system, we arrive at the integrability conditions in the form%
\begin{equation}
\begin{array}{c}
G_{11}u+G_{12}v=0, \\
G_{21}u+G_{22}v=0,%
\end{array}
\label{system with Gij}
\end{equation}%
where $G_{11},G_{12},G_{21}\;$and $G_{22}\;$are defined by formulas in
Theorem \ref{rank2}.

But for a $4$-web to be$\;$of rank two, one needs two independent solutions $%
u$ and $v.$ This proves (\ref{cond for rank 2}).
\end{proof}

\subsubsection{Examples}

\begin{example}
\label{rank2nonlin1}We consider the planar $4$-web formed by the coordinate
lines $y=\limfunc{const}.,\;x=\limfunc{const}.,\;$and by the level sets of
the functions%
\begin{equation*}
f(x,y)=x+y\;\text{and\ }g(x,y)=x^{2}+y^{2}.
\end{equation*}%
$($see Example $22)$.
\end{example}

We have already established that this $4$-web admits the abelian relation$%
\mathrm{\;}$
\begin{equation*}
u\omega _{1}+u\omega _{2}+u\omega _{3}=0,
\end{equation*}%
and its rank is either $1\;$or $2.$

In this case
\begin{eqnarray*}
c_{0} &=&-\frac{3(x-y)^{3}(x+y)}{xy^{5}}\neq 0,\; \\
\;c_{1} &=&\frac{x^{2}-y^{2}}{4x^{2}y^{3}},\ \;c_{2}=0; \\
\;c_{0,1} &=&-\frac{1}{2x^{3}},\;\ c_{0,2}=\frac{1}{2y^{3}}; \\
c_{1,1} &=&-\frac{1}{2x^{3}y},\;\ c_{1,2}=\frac{3x^{2}-y^{2}}{4x^{2}y^{4}};\;
\\
c_{2,1} &=&c_{2,2}=0,
\end{eqnarray*}%
and
\begin{equation*}
\;G_{11}=G_{12}=G_{21}=G_{22}=0.
\end{equation*}

It follows that conditions (\ref{cond for rank 2}) are satisfied for this $4$%
-web. Thus, the $4$-web is of rank two.

Two abelian relations for this web are:%
\begin{eqnarray*}
f_{1}+f_{2}-f_{3} &=&0, \\
f_{1}^{2}+f_{2}^{2}-f_{4} &=&0,
\end{eqnarray*}%
where%
\begin{equation*}
f_{1}=x,\;f_{2}=y,\;f_{3}=x+y,\text{\ }f_{4}=x^{2}+y^{2}.
\end{equation*}

\begin{example}
\label{rank2nonlin2}We consider the planar $4$-web formed by the coordinate
lines $y=\limfunc{const}.,\;x=\limfunc{const}.,\;$and by the level sets of
the functions%
\begin{equation*}
f(x,y)=\frac{x}{y}\;\text{and\ }g(x,y)=xy(x+y).
\end{equation*}
\end{example}

We have again $K=0,$ and one can check that the $4$-web linearizability
conditions \cite{AGL} are not satisfied. Therefore, our $4$-web is not
linearizable. By Theorem \ref{Poincare}$,\;$this $4$-web is not of maximum
rank three. Thus, the rank of the $4$-web in question can be $1$ or $2$.

We have%
\begin{eqnarray*}
c_{0} &=&\frac{3y^{3}(x^{2}-y^{2})}{2x(2x+y)^{2}(x+2y)^{2}}\neq 0;\; \\
\;c_{1} &=&\frac{3y^{5}(y-x)}{2x(2x+y)^{2}(x+2y)^{2}}\neq 0;\;c_{2}=0; \\
c_{0,1} &=&c_{0,2}=\frac{3y^{4}(2x^{4}-5x^{3}y-12x^{2}y^{2}-5xy^{3}+2y^{4})}{%
2x(2x+y)^{2}(x+2y)^{2}}; \\
c_{1,1} &=&c_{1,2}=-\frac{3y^{6}(4x^{3}-10x^{2}y-7xy^{2}+4y^{3})}{%
2x^{2}(2x+y)^{3}(x+2y)^{4}}; \\
\;c_{2,1} &=&c_{2,2}=0,
\end{eqnarray*}%
and as a result
\begin{equation*}
\;G_{11}=G_{12}=G_{21}=G_{22}=0.
\end{equation*}

Thus, the $4$-web is of rank two.

The two abelian relations are
\begin{equation*}
\renewcommand{\arraystretch}{1.5}%
\begin{array}{l}
\ln f_{1}-\ln f_{2}-\ln f_{3}=0, \\
\ln f_{1}+2\ln f_{2}+\ln \left( 1+f_{3}\right) -\ln f_{4}=0,%
\end{array}%
\renewcommand{\arraystretch}{1}
\end{equation*}%
where
\begin{equation*}
f_{1}=x,\;f_{2}=y,\;f_{3}=\frac{x}{y},\text{\ }f_{4}=xy(x+y).
\end{equation*}

Examples \ref{rank2nonlin1} and\ \ref{rank2nonlin2} lead us to the important
observation:

\begin{proposition}
\label{nonlin of webs with 2 abeq}In general$\;4$-webs of rank two are not
linearizable.
\end{proposition}

\subsection{$\mathbf{4}$-webs of rank one}

As we have seen before,$\;$a $4$-web can be of rank one if $c_{0}=0\;$but
one of the coefficients $c_{1}\;$and $c_{2}\;$of (\ref{eq4}) is not $0\;$or
if $c_{0}\neq 0.\;$The following theorem outlines the four cases when a $4$%
-web$\;$can be of rank one.

\begin{theorem}
\label{rank1}A planar $4$-web$\;$is of rank one if and only if one of the
following conditions holds:

\begin{enumerate}
\item $c_{0}=0,$ $J_{1}=J_{2}=0,$ where%
\begin{eqnarray*}
J_{1} &=&a_{2}c_{1}c_{2}(c_{1}-c_{2})+ac_{2}^{2}(c_{1,2}-c_{1,1}) \\
&&+c_{1}c_{2}(c_{1,1}+a(c_{2,1}-c_{1,2}-c_{2,2}))+c_{1}^{2}(ac_{2,2}-c_{2,1}),
\\
J_{2} &=&c_{1}^{2}\left( c_{1}-c_{2}\right) ^{2}K+\left(
c_{1,11}-c_{1,12}\right) c_{1}c_{2}\left( c_{2}-c_{1}\right) \\
&&+c_{1}^{2}\left( c_{1}-c_{2}\right) \left( c_{2,11}-c_{2,12}\right)
-c_{2}\left( 2c_{1}-c_{2}\right) c_{1,1}(c_{1,2}-c_{1,1}) \\
&&+c_{1}^{2}c_{2,1}(c_{1,2}-c_{2,2}+c_{2,1})+c_{1}^{2}c_{1,1}(c_{2,2}-2c_{2,1})
\end{eqnarray*}%
and $c_{1}\neq c_{2},\;c_{1}\neq 0.$

\item $c_{0}=0,c_{1}=c_{2}\neq 0,$ and $J_{3}=0,$ where%
\begin{equation*}
J_{3}=\left( a_{22}-a_{12}\right) \left( 1-a\right)
+a_{2}(a_{2}-a_{1})-\left( 1-a\right) ^{2}K.
\end{equation*}

\item $c_{0}=0,c_{1}=0,c_{2}\neq 0,$ and $J_{4}=0,$ where
\begin{equation*}
J_{4}=a_{12}a-a_{1}a_{2}-Ka^{2}.
\end{equation*}

\item $c_{0}\neq 0,$ and $J_{10}=J_{11}=J_{12}=0,$ where
\begin{equation*}
J_{10}=G_{11}G_{22}-G_{21}G_{12},
\end{equation*}%
\begin{eqnarray*}
J_{11} &=&c_{0}(G_{21,1}G_{22}-G_{22,1}G_{21})+(a_{2}c_{0}-ac_{1})G_{21}^{2}
\\
&&+(ac_{2}-a_{2}c_{0}+ac_{1})G_{21}G_{22}-ac_{2}G_{22}^{2},
\end{eqnarray*}%
\begin{eqnarray*}
J_{12} &=&c_{0}(G_{21,2}G_{22}-G_{22,2}G_{21})+(a_{2}c_{0}-ac_{1})G_{21}^{2}
\\
&&+a(c_{2}-c_{1})G_{21}G_{22}-c_{2}G_{22}^{2}.
\end{eqnarray*}
\end{enumerate}
\end{theorem}

\begin{proof}
First, we consider the case when $c_{0}=0,\;$ one of the coefficients $%
c_{1}\;$and $c_{2}\;$of (\ref{eq4}) is not $0,\;$and $4$-web $W_{4}\;$is of
rank one.

Then it follows from equation (\ref{eq4}) that
\begin{equation}
u=c_{1}t,\;\;v=-c_{2}t  \label{parametrization}
\end{equation}%
for some function $t.$

Differentiating these equations, we find that
\begin{equation*}
\begin{array}{ll}
u_{1}=c_{1,1}t+c_{1}t_{1}, & u_{2}=c_{1,2}t+c_{1}t_{2}, \\
v_{1}=-c_{2,1}t-c_{2}t_{1}, & v_{1}=-c_{2,2}t-c_{2}t_{2}.%
\end{array}%
\end{equation*}%
Substituting these expressions into the abelian equation, we get%
\begin{equation}
\begin{array}{l}
(c_{1,1}-c_{2,1})t+(c_{1}-c_{2})t_{1}=0, \\
(c_{1,1}-c_{1,2})t+c_{1}(t_{1}-t_{2})=0, \\
(-c_{2,1}+ac_{2,2}+a_{2}c_{2})t-c_{2}t_{1}+ac_{2}t_{2}=0.%
\end{array}
\label{t and ti}
\end{equation}%
If $c_{1}-c_{2}\neq 0\;$and\ $c_{1}\neq 0,\;$then solving the first two
equations with respect to $t_{1}\;$and $t_{2},$ we obtain%
\begin{equation}
t_{1}=\frac{t(c_{2,1}-c_{1,1})}{c_{1}-c_{2}},\;t_{2}=\frac{t(c_{1,1}-c_{2,1})%
}{c_{1}}+\frac{t(c_{2,1}-c_{1,1})}{c_{1}-c_{2}}.\;  \label{ti}
\end{equation}

Substituting these values into the last equation of the previous system, we
arrive at the equation $J_{1}=0,\;$where $J_{1}\;$is expressed as in Theorem %
\ref{rank1}.

Next, differentiating the third equation of (\ref{t and ti}) in the
direction $\left\{ \omega _{2}=0\right\} \;$and using symmetric derivatives,$%
\;$we find that%
\begin{eqnarray*}
&&\frac{t_{1}}{c_{1}}(c_{1,1}-c_{1,2})-\frac{t}{c_{1}^{2}}%
[(c_{1,11}-c_{1,12}+\frac{3Kc_{1}}{2})c_{1}-(c_{1,1}-c_{1,2})c_{1,1}] \\
&&-\frac{t_{1}}{c_{1}-c_{2}}(c_{1,1}-c_{2,1})-\frac{t}{(c_{1}-c_{2})^{2}}%
[(c_{1,11}-c_{2,11})(c_{1}-c_{2})-(c_{1,1}-c_{2,1}) \\
&&+\frac{t_{2}(c_{1,1}-c_{2,1})}{c_{1}-c_{2}}+t\frac{c_{1,12}-c_{2,12}}{%
(c_{1}-c_{2})^{2}}+t\frac{3K}{2(c_{1}-c_{2})}-t\frac{%
(c_{1,1}-c_{2,1})(c_{1,2}-c_{2,2})}{(c_{1}-c_{2})^{2}} \\
&=&0.
\end{eqnarray*}%
Substituting the values of $t_{1}\;$and $t_{2}\;$from (\ref{ti}) into the
above equation, we arrive at the equation $J_{2}=0\ $of Theorem \ref{rank1}.

Consider now the case: $c_{0}=0,\ c_{1}=c_{2}\neq 0.$

Then $u=-v,\;$and%
\begin{equation*}
\begin{array}{l}
u_{1}-u_{2}=0, \\
u_{1}-au_{2}-a_{2}u=0.%
\end{array}%
\end{equation*}%
Solving this system with respect to $u_{1}\;$and$\;u_{2},\;$we find that%
\begin{equation*}
u_{1}=\frac{a_{2}u}{1-a},\;\;u_{2}=\frac{a_{2}u}{1-a}.
\end{equation*}%
The compatibility of the above equations gives\ $J_{3}=0,$ where%
\begin{equation*}
J_{3}=\left( \frac{a_{2}u}{1-a}\right) _{2}-\left( \frac{a_{2}u}{1-a}\right)
_{1}-K=\frac{(a_{12}-a_{22})(a-1)-a_{2}(a_{1}-a_{2})}{(a-1)^{2}}-K.
\end{equation*}

Consider now the second excluded case: $c_{0}=0,\;c_{1}=0,\;c_{2}\neq 0.\;$

Then $u=0.\;$This implies
\begin{eqnarray*}
&&v_{1}=0, \\
&&av_{2}+a_{2}v=0,
\end{eqnarray*}%
and the compatibility condition $J_{4}=0,\;$where $%
J_{4}=a_{12}a-a_{1}a_{2}-Ka^{2}.$

Suppose now that $c_{0}\neq 0,\;$and a $4$-web$\;$is of rank one. Then the
abelian equation together with the compatibility condition $\varkappa =0$
gives the system%
\begin{eqnarray*}
u_{1} &=&-v_{1}, \\
u_{2} &=&-v_{1}, \\
u_{1} &=&-a_{2}v+\frac{1}{c_{0}}(c_{1}v+c_{2}u), \\
v_{2} &=&-\frac{1}{c_{0}}(c_{1}v+c_{2}u).
\end{eqnarray*}

Because the $4$-web$\;$is of rank one, system (\ref{system with Gij})\ has a
nonzero solution. Thus, its determinant vanishes:%
\begin{equation*}
J_{10}=G_{11}G_{22}-G_{12}G_{21}=0.
\end{equation*}%
Take, for example, the second equation of system (\ref{system with Gij}) and
differentiate it. Adding the resulting equations to the above system, we get
the conditions $J_{11}=J_{12}=0.$
\end{proof}

\subsubsection{Example}

\begin{example}
\label{rank1nonlin1}We consider the planar $4$-web$\;$formed by the
coordinate lines $y=\limfunc{const}.,\;x=\limfunc{const}t.$ and by the level
sets of the functions%
\begin{equation*}
f(x,y)=\frac{xy^{2}}{(x-y)^{2}}\;\text{and\ }g(x,y)=\frac{x^{2}y}{(x-y)^{2}}.
\end{equation*}
\end{example}

In this case, $c_{0}=0$ and $\;J_{1}=J_{2}=0.\;$Thus, we have the web of
type $1\;$as indicated in Theorem \ref{rank1},\ and this $4$-web is of rank
one.

The only abelian relation is
\begin{equation*}
\ln f_{1}-\ln f_{2}+\ln f_{3}-\ln f_{4}=0,
\end{equation*}%
where%
\begin{equation*}
f_{1}=x,\;f_{2}=y,\;f_{3}=\frac{xy^{2}}{(x-y)^{2}},\ f_{4}=\frac{x^{2}y}{%
(x-y)^{2}}.
\end{equation*}%
This example illustrates the following feature:

\begin{proposition}
\label{nonlinear webs of rank 1}In general, $4$-webs of rank one are not
linearizable.
\end{proposition}

\subsection{$\mathbf{4}$-webs with a constant basic invariant}

Consider first $4$-webs of maximum rank for which the basic invariant is
constant on one of web foliations.

Without loss of generality, we assume that $a$ is constant on the second
foliation, i.e.,
\begin{equation*}
a_{1}=0.
\end{equation*}%
Then, solving the system $c_{0}=c_{1}=c_{2}=0,$ we get
\begin{eqnarray*}
K &=&\frac{a_{2}^{2}}{4\left( a-1\right) ^{2}}-\frac{a_{22}}{4\left(
a-1\right) }, \\
K_{1} &=&\frac{a_{2}a_{22}}{2\left( a-1\right) ^{2}}-\frac{a_{2}^{3}}{%
2\left( a-1\right) ^{3}}, \\
K_{2} &=&\frac{a_{2}a_{22}}{4\left( a-1\right) ^{2}}-\frac{a_{2}^{3}}{%
4\left( a-1\right) ^{3}}.
\end{eqnarray*}

Differentiating the first equation and taking into account the remaining two
equations, we arrive at the conditions%
\begin{equation*}
a_{22}=\frac{a_{2}^{2}}{a-1},\;a_{222}=\frac{a_{2}^{3}}{\left( a-1\right)
^{2}}.
\end{equation*}%
The second equation above is the covariant derivative of the first one. The
first equation implies $K=0,$ and the curvatures of all other $3$-subwebs
vanish too. In other words, the $4$-web is a Mayrhofer web.

On the other hand, if we assume $a_{1}=0\ $and $K=0,$ then the only
condition for maximum rank is%
\begin{equation*}
a_{22}-\frac{a_{2}^{2}}{a-1}=0
\end{equation*}%
or%
\begin{equation*}
\delta _{2}\left( \frac{a_{2}}{a-1}\right) =0.
\end{equation*}

\begin{theorem}
\begin{enumerate}
\item A planar $4$-web with the basic invariant constant on one of the web
foliations is of maximum rank if one of its $3$-subwebs is parallelizable.
Then all other $3$-subwebs are parallelizable, and the $4$-web is a
Mayrhofer web.

\item If, say, $a_{1}=0$ and $K[1,2,3]=0$, then the $4$-web is of maximum
rank if $a_{2}/\left( a-1\right) $ is constant on the first foliation, i.e.,
if%
\begin{equation*}
\delta _{2}\left( \frac{a_{2}}{a-1}\right) =0.
\end{equation*}
\end{enumerate}
\end{theorem}

We conclude our discussion of $4$-webs of different ranks by the following
theorem whose proof utilizes all previous results.

\begin{theorem}
A planar $4$-web with a nonparallelizable $3$-subweb and a constant basic
invariant is of rank $0.$
\end{theorem}

\begin{proof}
Consider a $4$-web with a nonparallelizable $3$-subweb $[1,2,3]$ and a
constant basic invariant $a.\;$Then%
\begin{equation*}
c_{0}=K\neq 0,\;\;c_{1}=\frac{K_{1}-K_{2}}{4(a-1)},\;c_{2}=\frac{K_{1}-aK_{2}%
}{4a(a-1)}.
\end{equation*}%
Therefore, the web is not of rank three.

Substituting these expressions into equations of Theorem \ref{rank2}, we
find \
\begin{eqnarray*}
G_{11} &=&\frac{%
5(K_{1}^{2}-(a+1)K_{1}K_{2}+aK_{2}^{2})-4K(4K^{2}(a-1)+K_{11}-2aK_{12}+a^{2}K_{22})%
}{16(a-1)}, \\
G_{12} &=&\frac{5(K_{1}-K_{2})^{2}-4K(K_{11}-2K_{12}+K_{22})}{16(a-1)}, \\
G_{21} &=&\frac{-5(K_{1}-aK_{2})^{2}+4K(K_{11}-a(a-1)K_{12}+a^{3}K_{22})}{%
16a(a-1)}, \\
G_{22} &=&-\frac{%
5(K_{1}^{2}-(a+1)K_{1}K_{2}+aK_{2}^{2})+4K(4K^{2}(a-1)-K_{11}+(1+a)K_{12}-aK_{22})%
}{16(a-1)}.
\end{eqnarray*}

Thus, in general, $G_{ij}\neq 0$.\ Therefore, $4$-webs$\;$with a constant
basic invariant are not of rank two.

To check whether our $4$-web is of rank $1$, we note that $c_{0}=K\neq 0\;$%
because the $3$-subweb $[1,2,3]$ is nonparallelizable. Therefore, if the web
admits one abelian equation, then it would belong to class $4\;$of Theorem %
\ref{rank1}.$\;$This will be the case if $J_{10}=J_{12}=J_{12}=0.$

$\;$But
\begin{eqnarray*}
J_{10} &=&G_{11}G_{22}-G_{12}G_{21} \\
&=&\frac{1}{64}%
K[64K^{5}-5K_{2}^{2}K_{11}+16aK^{3}K_{22}-5(a+1)K_{1}^{2}K_{22}+4(a+1)KK_{11}K_{22}
\\
&&+5K_{1}K_{2}(K_{12}+aK_{22})-K_{12}(-5K_{1}^{2}+4K(4K^{2}+K_{11}+aK_{22}))].
\end{eqnarray*}

Thus in general, $J_{10}\neq 0.$\ As a result, the web does not admit even
one abelian equation, and it is of rank $0$.
\end{proof}

\section{Planar $\mathbf{5}$-webs$\;$}

\subsection{$\mathbf{5}$-webs of maximum rank}

Let us consider a planar $5$-web in the standard normalization%
\begin{equation*}
\begin{array}{l}
\omega _{1}+\omega _{2}+\omega _{3} =0, \\
a\omega _{1}+\omega _{2}+\;\omega_{4} =0, \\
b\omega _{1}+\omega _{2}+\;\omega_{5} =0,%
\end{array}%
\end{equation*}%
where $a$ and $b\;$are the basic invariants of the web.

The abelian relation for such a web has the form
\begin{equation*}
(w+au+bv)\omega _{1}+(w+u+v)\omega _{2}+w\omega _{3}+u\omega _{4}+v\omega
_{5}=0,
\end{equation*}%
where we have $\lambda _{1}=w+au+bv,\lambda _{2}=w+u+v,\lambda
_{3}=w,\lambda _{4}=u,$and $\ \lambda _{5}=v.$

The functions $w,u,$\ and$\ v$ satisfy the abelian equation%
\begin{equation*}
\begin{array}{l}
\delta _{1}\left( w\right) -\delta _{2}\left( w\right) =0, \\
\delta _{1}\left( u\right) -\delta _{2}\left( au\right) =0, \\
\delta _{1}\left( v\right) -\delta _{2}\left( bv\right) =0, \\
\delta _{1}\left( w\right) +\delta _{1}\left( u\right) +\delta _{1}\left(
v\right) =0,%
\end{array}%
\end{equation*}%
and their compatibility condition takes the form%
\begin{equation*}
\begin{array}{l}
\varkappa =\left( \Delta _{1}\Delta _{2}\Delta _{3}\delta _{1}-\delta
_{1}\Delta _{2}\Delta _{3}\Delta _{1}\right) \left( w\right) +\left( \Delta
_{1}\Delta _{2}\Delta _{3}\delta _{1}-\Delta _{1}\delta _{1}\Delta
_{3}\Delta _{2}\right) \left( u\right) \\
+\left( \Delta _{1}\Delta _{2}\Delta _{3}\delta _{1}-\Delta _{1}\Delta
_{2}\delta _{1}\Delta _{3}\right) \left( v\right) =0.%
\end{array}%
\end{equation*}

In the canonical coordinates in the jet bundles, the abelian equation has
the form
\begin{equation*}
\begin{array}{l}
w_{1}-w_{2}=0, \\
u_{1}-au_{2}-a_{2}u=0, \\
v_{1}-bv_{2}-b_{2}v=0, \\
u_{1}+v_{1}+w_{1}=0,%
\end{array}%
\end{equation*}%
and the obstruction $\varkappa $ equals%
\begin{equation*}
c_{0}w_{22}++c_{1}w_{2}+c_{2}v_{2}+c_{3}w+c_{4}u+c_{5}v=0,
\end{equation*}%
where
\begin{eqnarray*}
c_{0} &=&K+R[a,b]+R[b,a] \\
&&-\frac{a-a^{2}+b-b^{2}-4ab+2a^{2}b+2ab^{2}}{10\left( -1+a\right) a\left(
a-b\right) ^{2}\left( -1+b\right) b}a_{1}b_{1} \\
&&+\frac{2a-a^{2}+2b-b^{2}-4ab+a^{2}b+ab^{2}}{10\left( -1+a\right) \left(
a-b\right) ^{2}\left( -1+b\right) }a_{2}b_{2},
\end{eqnarray*}%
and%
\begin{eqnarray*}
R[a,b] &=&=\frac{\left( 1-3a+b\right) a_{11}+\left( 4a-3a^{2}-3b+4ab\right)
a_{12}+\left( a^{2}-3ab+a^{2}b\right) a_{22}}{10\left( -1+a\right) a\left(
a-b\right) } \\
&&+\frac{2a-6a^{2}+6a^{3}-b+4ab-6a^{2}b-b^{2}+2ab^{2}}{10\left( -1+a\right)
^{2}a^{2}\left( a-b\right) ^{2}}a_{1}^{2} \\
&&+\frac{4a^{2}-8a^{3}+3a^{4}-6ab+14a^{2}b-8a^{3}b+3b^{2}-6ab^{2}+4a^{2}b^{2}%
}{10\left( -1+a\right) ^{2}a^{2}\left( a-b\right) ^{2}}a_{1}a_{2} \\
&&+\frac{-a^{4}-2a^{2}b+6a^{3}b-a^{4}b-2a^{2}b^{2}}{10\left( -1+a\right)
^{2}a^{2}\left( a-b\right) ^{2}}a_{2}^{2}+\frac{-1+ab}{10\left( -1+a\right)
\left( a-b\right) ^{2}\left( -1+b\right) }a_{1}b_{2}
\end{eqnarray*}%
and the expressions for $c_{1},c_{2},c_{3},c_{4},$ and $c_{5}$ are given in
Section 7.2.

\begin{theorem}
A planar $5$-web is of maximum rank six if and only if the invariants $%
c_{0},c_{1},c_{2},c_{3},c_{4}\ $and $c_{5}$ vanish.
\end{theorem}

Note that $c_{0}$ contains the curvature function $K,$ while $c_{1},c_{2}$
contain $K$ and the linear combinations of covariant derivatives $%
K_{1},K_{2},$ and $c_{3},c_{4},c_{5}$ contain $K$ \ and the linear
combinations of the covariant derivatives $K_{1},K_{2}$ and the second
symmetrized covariant derivatives $K_{11},K_{12},K_{22}.$

\begin{corollary}
For $5$-webs of maximum rank, the curvature function $K$ of the $3$-subweb $%
[1,2,3]$ has the following expression in terms of the basic invariants $a$
and $b$ and their covariant derivatives:%
\begin{eqnarray*}
K &=&\frac{a-a^{2}+b-b^{2}-4ab+2a^{2}b+2ab^{2}}{10\left( -1+a\right) a\left(
a-b\right) ^{2}\left( -1+b\right) b}a_{1}b_{1} \\
&&-\frac{2a-a^{2}+2b-b^{2}-4ab+a^{2}b+ab^{2}}{10\left( -1+a\right) \left(
a-b\right) ^{2}\left( -1+b\right) }a_{2}b_{2}-R[a,b]-R[b,a].
\end{eqnarray*}
\end{corollary}

\subsection{Curvature of planar $\mathbf{5}$-webs}

Similar to the case of $4$-webs, the coefficient $c_{0}$ in the expression
of $\varkappa $ for $5$-webs has also an intrinsic geometric meaning,
Namely, let us define a curvature of$\ $a $5$-web as the arithmetic mean of
the curvatures of its $3$-subwebs. \ A planar $5$-web has ten $3$-subwebs:

$[1,2,3]$, $[1,2,4]$, $[1,3,4]$, $[2,3,4]$, $[1,2,5]$, $[1,3,5]$, $[2,3,5]$,
$[1,4,5]$, $[2,4,5]$, $[3,4,5]$.\

Normalizations with the curvature functions for $3$-subwebs $[1,2,3],[1,2,4]$%
,\newline
$[1,3,4]\;$and $[2,3,4]$ are given in Section 5.1. The similar expressions
for\newline
$3$-subwebs $[1,2,5],\;[1,3,5]$ and $[2,3,5]$ can be obtained by the
substitutions $a\rightarrow b$ and $4\rightarrow 5$:

\begin{itemize}
\item[{$[1,2,5]$:}] $\widetilde{\rho }_{1}=b\omega _{1},\widetilde{\rho }%
_{2}=\omega _{2},\widetilde{\rho }_{3}=\omega _{5},\widetilde{\rho }%
_{4}=\omega _{3}$;

\item[{$[1,3,5]$:}] $\widetilde{\sigma }_{1}=(b-1)\omega _{1},\widetilde{%
\sigma }_{2}=-\omega _{3},\widetilde{\sigma }_{3}=\omega _{5},\widetilde{%
\sigma }_{4}=-\omega _{2}$;

\item[{$[2,3,5]$:}] $\widetilde{\tau }_{1}=(b-1)\omega _{2},\widetilde{\tau }%
_{2}=b\omega _{3},\widetilde{\tau }_{3}=-\omega _{5},\widetilde{\tau }%
_{4}=b\omega _{1}$;
\end{itemize}

For the last three cases, we have

\begin{itemize}
\item[{$[1,4,5]$:}] $\zeta _{1}=(a-b)\omega _{1},\;\;\zeta _{2}=-a\omega
_{1}-\omega _{2},\;\zeta _{3}=b\omega _{1}+\omega _{2};$

\item[{$[2,4,5]$:}] $\eta _{1}=\frac{b-a}{b}\omega _{2},\;\;\eta
_{2}=-a\omega _{1}-\omega _{2},\;\ \;\eta _{3}=\frac{a}{b}(b\omega
_{1}+\omega _{2});$

\item[{$[3,4,5]$:}] $\theta _{1}=(a-b)\omega _{3},\;\theta _{2}=(1-b)(a\omega
_{1}+\omega _{2}),\;\;\theta _{3}=(a-1)(b\omega _{1}+\omega _{2});$
\end{itemize}

For the curvature functions $K[l,m,n]$,\ we get the following expressions:
\begin{eqnarray*}
K[1,2,5] &=&\frac{1}{b}\left( K-\frac{b_{12}}{b}+\frac{b_{1}b_{2}}{b^{2}}%
\right) , \\
K[1,3,5] &=&\frac{1}{b-1}\left[ K+\frac{b_{2}(b_{1}-b_{2})}{(1-b)^{2}}+\frac{%
b_{12}-b_{22}}{1-b}\right] , \\
K[2,3,5] &=&\frac{1}{b(b-1)}\left[ K+\frac{(2b-1)b_{1}(b_{1}-b_{2})}{%
b^{2}(1-b)^{2}}+\frac{b_{11}-b_{12}}{b(1-b)}\right]
\end{eqnarray*}%
and%
\begin{eqnarray*}
K[1,4,5] &=&\frac{K-a_{22}}{b-a}+\frac{b_{12}-a_{12}+a(a_{22}-b_{22})}{%
\left( a-b\right) ^{2}} \\
&&+\frac{(b_{2}-a_{2})(a_{2}b-ab_{2}-a_{1}+b_{1})}{(a-b)^{3}},
\end{eqnarray*}%
\begin{eqnarray*}
K[2,4,5] &=&\frac{bK}{a(b-a)}-\frac{aa_{12}-a_{1}a_{2}}{a^{3}(b-a)}+\frac{%
a_{11}b-ab_{11}}{a^{2}b(b-a)^{2}} \\
&&-\frac{a_{12}b+a_{1}b_{2}-a_{2}b_{1}-ab_{12}}{ab(b-a)^{2}}+\frac{%
(a_{1}b-ab_{1})(2bb_{2}-ab_{2}-a_{2}b)}{ab^{2}(b-a)^{3}} \\
&&-\frac{(a_{1}b-ab_{1})(a_{1}b^{2}+2ab(b_{1}-a_{1})-a^{2}b_{1})}{%
a^{3}b^{2}(b-a)^{3}},
\end{eqnarray*}%
\begin{eqnarray*}
K[3,4,5] &=&\frac{K}{(a-1)(b-1)(b-a)}-\frac{aa_{2}(b_{1}-a_{1}-b_{2}+a_{2})}{%
(a-1)^{3}(b-1)(b-a)^{2}} \\
&&+\frac{a_{2}(ab_{2}-b_{1}+(b-1)a_{2})}{(a-1)^{3}(b-1)^{2}(b-a)}+\frac{%
a_{1}\left( b_{1}-a_{1}-b_{2}+a_{2}\right) }{(a-1)^{3}(b-1)(b-a)^{2}} \\
&&+\frac{a_{1}\left( ab_{2}-b_{1}+(b-1)a_{2}\right) }{(a-1)^{3}(b-1)^{2}(b-a)%
}-\frac{a(b_{12}-a_{12}-b_{22}+a_{22})}{(a-1)^{2}(b-1)(b-a)^{2}} \\
&&-\frac{a_{2}(b_{1}-a_{1}-b_{2}+a_{2}}{(a-1)^{2}(b-1)(b-a)^{2}}+\frac{%
a(b_{1}-a_{1}-b_{2}+a_{2})(b_{2}-a_{2})}{(a-1)^{2}(b-1)(b-a)^{3}} \\
&&-\frac{(b_{1}-ab_{2})b_{2}}{(a-1)^{2}(b-1)^{3}(b-a)}-\frac{%
a_{2}b_{2}-b_{12}+ab_{22}+(b-1)a_{22}}{(a-1)^{2}(b-1)^{2}(b-a)} \\
&&+\frac{(b_{1}-a_{1}-b_{2}+a_{2})(b_{1}-a_{1})}{(a-1)^{2}(b-1)(b-a)^{3}}-%
\frac{b_{11}-a_{11}-b_{12}+a_{12}}{(a-1)^{2}(b-1)(b-a)^{2}} \\
&&-\frac{-b_{11}+a_{1}b_{2}+ab_{12}+(b-1)a_{12}}{(a-1)^{2}(b-1)^{2}(b-a)}-%
\frac{(b_{1}-ab_{2})b_{1}}{(a-1)^{2}(b-1)^{3}(b-a)},
\end{eqnarray*}

Define a \emph{curvature }$2$\emph{-form }$L\omega _{1}\wedge \omega _{2}$%
\emph{\ of the }$5$\emph{-web} as follows%
\begin{eqnarray*}
10L\omega _{1}\wedge \omega _{2} &=&K[1,2,3]\omega _{1}\wedge \omega _{2} \\
&&+K[1,2,4]\rho _{1}\wedge \rho _{2}+K[1,3,4]\sigma _{1}\wedge \sigma
_{2}+K[2,3,4]\tau _{1}\wedge \tau _{2} \\
&&+K[1,2,5]\widetilde{\rho }_{1}\wedge \widetilde{\rho }_{2}+K[1,3,5]%
\widetilde{\sigma }_{1}\wedge \widetilde{\sigma }_{2}+K[2,3,5]\widetilde{%
\tau }_{1}\wedge \widetilde{\tau }_{2} \\
&&+K[1,4,5]\zeta _{1}\wedge \zeta _{2}+K[2,4,5]\eta _{1}\wedge \eta
_{2}+K[3,4,5]\theta _{1}\wedge \theta _{2}
\end{eqnarray*}

The straightforward calculation shows that
\begin{eqnarray*}
10L &=&K+aK[1,2,4]+(a-1)K[1,3,4] \\
&&+a(a-1)K[2,3,4]+bK[1,2,5]+(b-1)K[1,3,5]+b(b-1)K[2,3,5] \\
&&+(b-a)K[1,4,5]+\frac{b}{a(b-a)}K[2,4,5]+(a-1)(b-1)(b-a)K[3,4,5]
\end{eqnarray*}%
and
\begin{equation*}
L=c_{0}.
\end{equation*}

\begin{theorem}
The curvature of a planar $5$-web of maximum rank equals zero.
\end{theorem}

\subsection{$\mathbf{5}$-webs with constant basic invariants}

We conclude this section by consideration of $5$-webs with constant basic
invariants. One can check that in this case $c_{0}$ coincides with the
curvature function $K$ of the $3$-subweb $[1,2,3],$ and the expressions for $%
c_{1},c_{2},c_{3},c_{4}\ $and$\ c_{5}$ are linear combinations of $%
K_{i},K_{ij}$ and $K^{2}.$ In other words, $%
c_{0}=c_{1}=c_{2}=c_{3}=c_{4}=c_{5}=0$ if and only if $K=0.$

\begin{theorem}
A planar $5$-web with constant basic invariants is of maximum rank if and
only if it is parallelizable.
\end{theorem}

\subsection{Examples}

In this section we consider two examples of $5$-webs.

\begin{example}
\label{Bol 5-web}We consider the Bol $5$-web formed by the coordinate lines $%
y=\limfunc{const}.,\;x=\limfunc{const}.,\;$and by the level sets of the
functions%
\begin{equation*}
f(x,y)=\frac{x}{y},\text{\ }g_{4}(x,y)=\frac{1-y}{1-x}\;\text{and }%
g_{5}(x,y)=\frac{x-xy}{y-xy}
\end{equation*}%
$($see Example $8$ in Section $5.2\;$of$\;\cite{AGL})$.\
\end{example}

First note that because the $3$-subweb $[1,2,3]$ of this $5$-web is
parallelizable, we have $K=0.$

For this $5$-web, we have%
\begin{equation*}
a=\frac{x(y-1)}{y(x-1)},b=\frac{y-1}{x-1},
\end{equation*}%
and
\begin{equation*}
c_{i}=0
\end{equation*}%
for $i=0,1,...,5.$

Thus, the Bol $5$-web is of maximum rank.

Using the linearizability conditions for planar $5$-webs given in \cite{AGL}%
, we observe that the Bol $5$-web is not linearizable. Indeed, one of the
invariants, namely, the second-order invariant
\begin{equation*}
\mu _{\lbrack 1,2,3,4]}-\mu _{\lbrack 1,2,3,5]}=\frac{\partial
_{1}a-a\partial _{2}a}{a-a^{2}}-\frac{\partial _{1}b-b\partial _{2}b}{b-b^{2}%
}
\end{equation*}%
does not vanish.

For this $5$-web, we have%
\begin{eqnarray*}
\omega _{1} &=&-\frac{1}{y}dx,\;\;\omega _{2}=\frac{x}{y^{2}}dy,\ \ \omega
_{3}=\frac{1}{y}dx-\frac{x}{y^{2}}dy, \\
\omega _{4} &=&\frac{x(y-1)}{(x-1)y^{2}}dx-\frac{x}{y^{2}}dy,\ \omega _{5}=%
\frac{x(y-1)}{(x-1)y}dx-\frac{x^{2}}{y^{2}}dy,
\end{eqnarray*}%
and solving the abelian equation, we get the following abelian relations:%
\begin{equation*}
\renewcommand{\arraystretch}{1.5}%
\begin{array}{l}
\ln f_{1}-\ln f_{2}-\ln f_{3}=0, \\
\ln f_{3}+\ln f_{4}-\ln f_{5}=0, \\
\ln (1-f_{1})-\ln \left( 1-f_{2}\right) +\ln f_{4}=0, \\
\ln (1-f_{1})-\ln \left( 1-f_{3}\right) +\ln (1-f_{5})=0, \\
\ln \frac{1-f_{1}}{f_{1}}-\ln \frac{1-f_{3}}{f_{3}}+\ln (1-f_{4})=0, \\
\mathbf{D}_{2}(f_{1})-\mathbf{D}_{2}(f_{2})-\mathbf{D}_{2}(f_{3})-\mathbf{D}%
_{2}(f_{4})+\mathbf{D}_{2}(f_{5})=0,%
\end{array}%
\renewcommand{\arraystretch}{1}
\end{equation*}%
where%
\begin{equation*}
f_{1}=x,\;f_{2}=y,\;f_{3}=\frac{x}{y},\text{\ }f_{4}=\frac{1-y}{1-x},\;\text{
}f_{5}=\frac{x-xy}{y-xy}
\end{equation*}%
and
\begin{equation*}
\mathbf{D}_{2}(t)=\mathbf{Li}_{2}t+\frac{1}{2}\ln t\;\ln (1-t)-\frac{\pi ^{2}%
}{6},\;\;0<t<1;\;\mathbf{Li}_{2}t=\sum_{n=1}^{\infty }\frac{t^{n}}{n^{2}},
\end{equation*}%
or%
\begin{equation*}
\mathbf{D}_{2}(t)=-\frac{1}{2}\int_{0}^{t}\left( \frac{\ln |1-s|}{s}+\frac{%
\ln |s|}{1-s}\right) ds-\frac{\pi ^{2}}{6},\;\;0<t<1,
\end{equation*}%
is the version of the original Rogers dilogarithm (see \cite{Rog 07} and %
\cite{Le 81})\ normalized so that RHS of the last abelian relation is$\;0$.

This example leads us to the following important observation

\begin{proposition}
In general, planar $5$-webs of maximum rank are not linearizable $($%
algebraizable$)$.
\end{proposition}

\begin{example}
We consider the $5$-web formed by the coordinate lines $y=\limfunc{const}%
.,\;x=\limfunc{const}.,\;$and by the level sets of the functions%
\begin{equation*}
f(x,y)=\frac{x}{y},\text{\ }g_{4}(x,y)=\frac{x}{(x-1)(y-1)}\;\text{and }%
g_{5}(x,y)=\frac{y}{(x-1)(y-1)}.
\end{equation*}%
This is the $5$-subweb $[1,2,3,7,8]\;$of the $29$-web $K(5)\;($see \cite{P
04a}, Section $7.2.3)$.\
\end{example}

As in Example \ref{Bol 5-web}, we have
\begin{equation*}
K=0,\ a=\frac{1-y}{y(x-1)},\ b=\frac{x(1-y)}{x-1},
\end{equation*}%
and
\begin{equation*}
c_{i}=0
\end{equation*}%
for $i=0,1,...,5.$

Thus, the planar $5$-web in question is of maximum rank.

One can also check the linearizability conditions (see of \cite{AGL}, p.
445) hold, and this $5$-web is linearizable (algebraizable). Therefore, the $%
5$-web in question is algebraizable $5$-web of maximum rank.

For this $5$-web, we have%
\begin{eqnarray*}
\omega _{1} &=&-\frac{1}{y}dx,\;\;\omega _{2}=\frac{x}{y^{2}}dy,\ \ \omega
_{3}=\frac{1}{y}dx-\frac{x}{y^{2}}dy, \\
\omega _{4} &=&\frac{1-y}{(x-1)y^{2}}dx-\frac{x}{y^{2}}dy,\ \omega _{5}=%
\frac{x(1-y)}{(x-1)y^{2}}dx-\frac{x}{y^{2}}dy,
\end{eqnarray*}%
and solving the abelian equation, we get the following abelian relations%
\begin{equation*}
\renewcommand{\arraystretch}{1.5}
\begin{array}{l}
\ln f_{1}-\ln f_{2}-\ln f_{3} =0, \\
\ln f_{3}-\ln f_{4}+\ln f_{5} =0, \\
\ln \frac{f_{1}}{1-f_{1}}-\ln \left( 1-f_{2}\right) -\ln f_{4} =0, \\
\ln (1-f_{1})-\ln \frac{f_{2}}{1-f_{2}}+\ln f_{5} =0, \\
\frac{f_{1}-1}{f_{1}}+f_{2}-f_{3}-\frac{1}{f_{4}} =0, \\
f_{1}-\frac{1-f_{2}}{f_{2}}+f_{3}+\frac{1}{f_{5}} =0,%
\end{array}%
\renewcommand{\arraystretch}{1}
\end{equation*}

where%
\begin{equation*}
f_{1}=x,\;f_{2}=y,\;f_{3}=\frac{x}{y},\text{\ }f_{4}=\frac{x}{(1-x)(1-y)},\;%
\text{ }f_{5}=\frac{y}{(1-x)(1-y)}.
\end{equation*}

\section{Appendix}

\subsection{Multi-brackets}

Compatibility conditions for linear (as well as nonlinear) PDEs systems can
be expressed in terms of multi-brackets introduced in \cite{KL 06}. Here we
give the necessary formulae for multi-brackets which we use in this paper.

Consider a linear PDEs system of $n+1$ differential equations for $n$
unknown functions:%
\begin{equation*}
\begin{Vmatrix}
a_{11} & \cdots & a_{1n} \\
\cdot & \cdots & \cdot \\
a_{n+11} & \cdots & a_{n+1~n}%
\end{Vmatrix}%
\begin{Vmatrix}
u_{1} \\
\vdots \\
u_{n}%
\end{Vmatrix}%
=0,
\end{equation*}%
where $a_{ij}$ are linear differential operators.

Then the multi-bracket $\left\{ a_{1},...,a_{n+1}\right\} $ of scalar
differential operators $a_{i}=\left( a_{i1},...,a_{in}\right) $ is given by
the formula%
\begin{equation*}
\left\{ a_{1},...,a_{n+1}\right\} =\dsum\limits_{i=1}^{n+1}\left( -1\right)
^{i-1}\limfunc{Ndet}\left( A_{i}\right) a_{i},
\end{equation*}%
where $A_{i}$ is the $n\times n$ matrix obtained from the matrix $\left\|
a_{ij}\right\| $ by deleting $i$-th row, and $\limfunc{Ndet}$ is a
non-commutative version of the determinant function produced by the standard
formula of \ decomposition of the determinant with respect to the first
column.

The abelian equation is the equation of this form and it satisfies the
conditions of \cite{KL 06}. Thus, the multi-brackets give the compatibility
conditions for the abelian equations.

\subsection{Coefficients for obstruction}

Below we list coefficients $c_{1},c_{2},c_{3},c_{4}\;$and$\;c_{5}.$ They are
\begin{equation*}
c_{1}=\frac{j_{1}}{3ad},\ c_{2}=\frac{j_{2}}{3ad},c_{3}=\frac{j_{3}}{d}%
,c_{4}=\frac{j_{4}}{6a},c_{5}=\frac{j_{5}}{6a},
\end{equation*}%
where
\begin{equation*}
d=-10a^{2}b^{2}\left( a-1\right) ^{2}\left( b-1\right) ^{2}\left( a-b\right)
^{2}
\end{equation*}%
and

\begin{equation*}
\begin{array}{ll}
j_{1}= & -3(-1+b)b^{2}(6a^{3}-6a^{2}(1+b)-b(1+b)+2a(1+b)^{2})a_{1}^{3} \\
& -3ba_{1}^{2}(-(-1+b)b(-6a^{4}-b^{2}(5+2b)+a^{3}(16+13b) \\
& +ab(13+12b+4b^{2})-a^{2}(10+23b+12b^{2}))a_{2}+a((b+a^{2}(7-12b) \\
& b-2b^{3}+a^{3}(-2+4b)+a(1-8b+7b^{2}+4b^{3}))b_{1}+b(-(-2+b) \\
& b^{2}+a^{3}(-6+4b)+a^{2}(3+7b-5b^{2})+a(1-6b-b^{2}+2b^{3})) b_{2})) \\
& -a_{1}(3(-1+b)b^{2}(-6b^{3}-3a^{4}(2+b)+ab^{2}(17+8b) \\
& -a^{2}b(16+23b)+a^{3}(4+24b+b^{2}))a_{2}^{2}+3(-1+a) \\
& aba_{2}((2(2-3b)b^{2}+a^{3}(-2+4b)+a^{2}(1+5b-12b^{2}) \\
& +ab(-8+9b+5b^{2}))b_{1}+b(-2(-2+b)b^{2}+a^{3}(-6+4b) \\
& +a^{2}(7+3b-4b^{2})+ab(-8-b+3b^{2}))b_{2}) \\
&
+a(-31a^{3}b^{2}K+28a^{4}b^{2}K+3a^{5}b^{2}K+50a^{2}b^{3}K+22a^{3}b^{3}K-69a^{4}b^{3}K
\\
& -3a^{5}b^{3}K-19ab^{4}K-92a^{2}b^{4}K+70a^{3}b^{4}K+41a^{4}b^{4}K+42ab^{5}K
\\
& +19a^{2}b^{5}K-61a^{3}b^{5}K-23ab^{6}K+23a^{2}b^{6}K-3(-1+b)b^{2}(6a^{3}
\\
& -b(1+2b)-3a^{2}(2+3b)+a(1+8b+3b^{2}))a_{11}+3(-1+b)b^{2}(-6a^{4} \\
& -b^{2}(5+2b)+3a^{3}(4+5b)+ab(13+12b+4b^{2})-a^{2}(8+22b+13b^{2}))a_{12}%
\end{array}%
\end{equation*}

\newpage

\begin{equation*}
\begin{array}{ll}
&
+6a^{3}b^{2}a_{22}-6a^{4}b^{2}a_{22}-24a^{2}b^{3}a_{22}+21a^{3}b^{3}a_{22}+27ab^{4}a_{22}-6a^{2}b^{4}a_{22}
\\
&
-18a^{3}b^{4}a_{22}+6a^{4}b^{4}a_{22}-9b^{5}a_{22}-18ab^{5}a_{22}+27a^{2}b^{5}a_{22}-9a^{3}b^{5}a_{22}
\\
&
+9b^{6}a_{22}-9ab^{6}a_{22}+3a^{2}b^{6}a_{22}+3a^{2}b_{1}^{2}-3a^{4}b_{1}^{2}-6abb_{1}^{2}-6a^{2}bb_{1}^{2}+6a^{3}bb_{1}^{2}
\\
&
+6a^{4}bb_{1}^{2}+18ab^{2}b_{1}^{2}-18a^{3}b^{2}b_{1}^{2}-18ab^{3}b_{1}^{2}+18a^{2}b^{3}b_{1}^{2}-9a^{3}b_{1}b_{2}+9a^{4}b_{1}b_{2}
\\
&
+18a^{2}bb_{1}b_{2}-18a^{4}bb_{1}b_{2}-12ab^{2}b_{1}b_{2}-30a^{2}b^{2}b_{1}b_{2}+30a^{3}b^{2}b_{1}b_{2}
\\
&
+12a^{4}b^{2}b_{1}b_{2}+24ab^{3}b_{1}b_{2}-24a^{3}b^{3}b_{1}b_{2}-9ab^{4}b_{1}b_{2}+9a^{2}b^{4}b_{1}b_{2}
\\
&
+6a^{2}b^{2}b_{2}^{2}-6a^{4}b^{2}b_{2}^{2}-18a^{2}b^{3}b_{2}^{2}+18a^{3}b^{3}b_{2}^{2}+3ab^{4}b_{2}^{2}-3a^{3}b^{4}b_{2}^{2}-3a^{2}bb_{11}
\\
&
+3a^{4}bb_{11}+3ab^{2}b_{11}+12a^{2}b^{2}b_{11}-12a^{3}b^{2}b_{11}-3a^{4}b^{2}b_{11}-12ab^{3}b_{11}
\\
&
+12a^{3}b^{3}b_{11}+9ab^{4}b_{11}-9a^{2}b^{4}b_{11}+9a^{3}bb_{12}-9a^{4}bb_{12}-21a^{2}b^{2}b_{12}
\\
&
+21a^{4}b^{2}b_{12}+12ab^{3}b_{12}+30a^{2}b^{3}b_{12}-30a^{3}b^{3}b_{12}-12a^{4}b^{3}b_{12}-21ab^{4}b_{12}
\\
&
+21a^{3}b^{4}b_{12}+9ab^{5}b_{12}-9a^{2}b^{5}b_{12}+9a^{3}b^{2}b_{22}-9a^{4}b^{2}b_{22}-12a^{2}b^{3}b_{22}
\\
&
+12a^{4}b^{3}b_{22}+3ab^{4}b_{22}+12a^{2}b^{4}b_{22}-12a^{3}b^{4}b_{22}-3a^{4}b^{4}b_{22}-3ab^{5}b_{22}
\\
&
+3a^{3}b^{5}b_{22}))-a(-3a^{2}(2+a(-3+b))(-1+b)b^{3}a_{2}^{3}+3a^{2}ba_{2}^{2}((a^{2}(1-2b)
\\
& +b+b^{2}-3b^{3}+ab(-3+4b+b^{2}))b_{1}+(-2+b)b(a^{2}(-2+b)+b-2b^{2} \\
& +a(1+b^{2}))b_{2})+.a_{2}(-3(-1+b)b^{2}(3a^{4}-7a^{3}(1+b)+2b^{2}(1+b) \\
& -ab(5+9b+2b^{2})+a^{2}(3+14b+6b^{2}))a_{11}+3(-1+b)b^{2}(6b^{3} \\
& +3a^{4}(2+b)-ab^{2}(17+8b)+a^{2}b(15+26b)-a^{3}(4+24b+3b^{2}))a_{12} \\
& +a(3a(-1+b)b^{3}(a^{2}(-4+b)-3b+a(3+4b-b^{2}))a_{22}+ \\
& (-1+a)(3a(a^{2}(-1+2b)+a(-1+4b-6b^{2})+2b(1-3b+3b^{2}))b_{1}^{2} \\
& +3a(b^{2}(4-8b+3b^{2})-2ab(3-7b+4b^{2})+a^{2}(3-6b+4b^{2}))b_{1}b_{2} \\
& +b(-3ab(2a^{2}+b^{2}+a(2-6b+b^{2}))b_{2}^{2}+(a-b)(-1+b) \\
& (-3a(1+a-3b)b_{11}-3a((4-3b)b+a(-3+4b))b_{12} \\
& -b((a^{2}(23-19b)+3b^{2}+ab(-38+31b))K+3a(a(-3+b)+b)b_{22})))))) \\
& +3(-1+a)a(a-b)b((a-b)(-1+b)ba_{111}+(-1+b)b(3a^{2}+2b(1+b) \\
&
-a(2+5b))a_{112}+2a^{2}ba_{122}-5ab^{2}a_{122}+3b^{3}a_{122}+3ab^{3}a_{122}-2a^{2}b^{3}a_{122}
\\
&
-3b^{4}a_{122}+2ab^{4}a_{122}-a^{2}b^{2}a_{222}+ab^{3}a_{222}+a^{2}b^{3}a_{222}-ab^{4}a_{222}+4a^{2}Kb_{1}
\\
&
-4a^{3}Kb_{1}-8abKb_{1}+8a^{3}bKb_{1}+12ab^{2}Kb_{1}-12a^{2}b^{2}Kb_{1}+aa_{11}b_{1}
\\
&
-2ba_{11}b_{1}-2aba_{11}b_{1}+3b^{2}a_{11}b_{1}-2aa_{12}b_{1}+2a^{2}a_{12}b_{1}+4ba_{12}b_{1}
\\
&
-4a^{2}ba_{12}b_{1}-6b^{2}a_{12}b_{1}+6ab^{2}a_{12}b_{1}-a^{2}a_{22}b_{1}+2aba_{22}b_{1}+2a^{2}ba_{22}b_{1}
\\
&
-3ab^{2}a_{22}b_{1}+12a^{2}bKb_{2}-12a^{3}bKb_{2}-8ab^{2}Kb_{2}+8a^{3}b^{2}Kb_{2}
\\
&
+4ab^{3}Kb_{2}-4a^{2}b^{3}Kb_{2}+3aba_{11}b_{2}-2b^{2}a_{11}b_{2}-2ab^{2}a_{11}b_{2}+b^{3}a_{11}b_{2}
\\
&
-6aba_{12}b_{2}+6a^{2}ba_{12}b_{2}+4b^{2}a_{12}b_{2}-4a^{2}b^{2}a_{12}b_{2}-2b^{3}a_{12}b_{2}+2ab^{3}a_{12}b_{2}
\\
&
-3a^{2}ba_{22}b_{2}+2ab^{2}a_{22}b_{2}+2a^{2}b^{2}a_{22}b_{2}-ab^{3}a_{22}b_{2}-5a^{2}bK_{1}+5a^{3}bK_{1}
\\
&
+5ab^{2}K_{1}-5a^{3}b^{2}K_{1}-5ab^{3}K_{1}+5a^{2}b^{3}K_{1}+5a^{2}b^{2}K_{2}-5a^{3}b^{2}K_{2}
\\
& -5ab^{3}K_{2}+5a^{3}b^{3}K_{2}+5ab^{4}K_{2}-5a^{2}b^{4}K_{2}));%
\end{array}%
\end{equation*}

\begin{equation*}
\begin{array}{ll}
j_{2}= & -3(-1+b)b^{3}(6a^{3}-6a^{2}(1+b)-b(1+b)+2a(1+b)^{2})a_{1}^{3} \\
& +3b^{2}a_{1}^{2}((-1+b)b(-6a^{4}-b^{2}(5+2b)+a^{3}(16+13b) \\
& +ab(13+12b+4b^{2})-a^{2}(10+23b+12b^{2}))a_{2}+a((2a^{3}(-2+b)+b^{3} \\
& +a^{2}(6-7b+6b^{2})-a(3-6b+5b^{2}+2b^{3}))b_{1}+(-1+a)b(2a^{2}b \\
& +b(-1+2b)-a(-3+5b+b^{2}))b_{2}))+ba_{1}(-3(-1+b)b^{2}(-6b^{3} \\
& -3a^{4}(2+b)+ab^{2}(17+8b)-a^{2}b(16+23b)+a^{3}(4+24b+b^{2}))a_{2}^{2} \\
& -3aba_{2}((a^{4}(1+b)+b^{2}(-1+3b)+ab(2-5b-5b^{2})+a^{3}(-5+b-4b^{2}) \\
&
+a^{2}(3-3b+11b^{2}+b^{3}))b_{1}-b(-a^{4}(-3+b)+b^{2}(1+b)+a^{3}(-5+b-4b^{2})%
\end{array}%
\end{equation*}

\newpage

\begin{equation*}
\begin{array}{ll}
&
-ab(2+5b+b^{2})+a^{2}(3+b+7b^{2}+b^{3}))b_{2})+a(7a^{3}b^{2}K+8a^{4}b^{2}K-3a^{5}b^{2}K
\\
&
-14a^{2}b^{3}K-34a^{3}b^{3}K-3a^{4}b^{3}K+3a^{5}b^{3}K+7ab^{4}K+44a^{2}b^{4}K+26a^{3}b^{4}K
\\
&
-5a^{4}b^{4}K-18ab^{5}K-31a^{2}b^{5}K+a^{3}b^{5}K+11ab^{6}K+a^{2}b^{6}K-3a^{2}(2+3b)
\\
& +3(-1+b)b^{2}(6a^{3}-b(1+2b)+a(1+8b+3b^{2}))a_{11} \\
& +3(-1+b)b^{2}(6a^{4}+b^{2}(5+2b)-3a^{3}(4+5b)-ab(13+12b+4b^{2}) \\
&
+a^{2}(8+22b+13b^{2}))a_{12}-6a^{3}b^{2}a_{22}+6a^{4}b^{2}a_{22}+24a^{2}b^{3}a_{22}-21a^{3}b^{3}a_{22}
\\
&
-27ab^{4}a_{22}+6a^{2}b^{4}a_{22}+18a^{3}b^{4}a_{22}-6a^{4}b^{4}a_{22}9b^{5}a_{22}+18ab^{5}a_{22}-27a^{2}b^{5}a_{22}
\\
&
+9a^{3}b^{5}a_{22}-9b^{6}a_{22}+9ab^{6}a_{22}-3a^{2}b^{6}a_{22}-3a^{4}b_{1}^{2}+9abb_{1}^{2}-18a^{2}bb_{1}^{2}
\\
&
+15a^{3}bb_{1}^{2}+6a^{4}bb_{1}^{2}-18ab^{2}b_{1}^{2}+21a^{2}b^{2}b_{1}^{2}-18a^{3}b^{2}b_{1}^{2}+12ab^{3}b_{1}^{2}-6a^{2}b^{3}b_{1}^{2}
\\
&
-3a^{3}b_{1}b_{2}+9a^{4}b_{1}b_{2}+6a^{2}bb_{1}b_{2}-15a^{3}bb_{1}b_{2}-15a^{4}bb_{1}b_{2}+9ab^{2}b_{1}b_{2}-9a^{2}b^{2}b_{1}b_{2}
\\
&
+33a^{3}b^{2}b_{1}b_{2}+3a^{4}b^{2}b_{1}b_{2}-15ab^{3}b_{1}b_{2}+3a^{2}b^{3}b_{1}b_{2}-12a^{3}b^{3}b_{1}b_{2}+3ab^{4}b_{1}b_{2}
\\
&
+3a^{2}b^{4}b_{1}b_{2}-3a^{2}b^{2}b_{2}^{2}+3a^{3}b^{2}b_{2}^{2}-3a^{4}b^{2}b_{2}^{2}+9a^{2}b^{3}b_{2}^{2}-6a^{3}b^{3}b_{2}^{2}+3a^{4}b^{3}b_{2}^{2}
\\
&
-6a^{2}b^{4}b_{2}^{2}+3a^{3}b^{4}b_{2}^{2}+3a^{2}bb_{11}-6a^{3}bb_{11}+6a^{4}bb_{11}-3ab^{2}b_{11}-3a^{2}b^{2}b_{11}
\\
&
+3a^{3}b^{2}b_{11}-6a^{4}b^{2}b_{11}+9ab^{3}b_{11}-3a^{2}b^{3}b_{11}+3a^{3}b^{3}b_{11}-6ab^{4}b_{11}+3a^{2}b^{4}b_{11}
\\
&
+3a^{3}bb_{12}-9a^{4}bb_{12}+3a^{2}b^{2}b_{12}+6a^{3}b^{2}b_{12}+15a^{4}b^{2}b_{12}-6ab^{3}b_{12}-6a^{2}b^{3}b_{12}
\\
&
-18a^{3}b^{3}b_{12}-6a^{4}b^{3}b_{12}+9ab^{4}b_{12}+6a^{2}b^{4}b_{12}+9a^{3}b^{4}b_{12}-3ab^{5}b_{12}-3a^{2}b^{5}b_{12}
\\
&
-3a^{3}b^{2}b_{22}+3a^{2}b^{3}b_{22}+9a^{3}b^{3}b_{22}-3a^{4}b^{3}b_{22}-9a^{2}b^{4}b_{22}-3a^{3}b^{4}b_{22}+3a^{4}b^{4}b_{22}
\\
&
+6a^{2}b^{5}b_{22}-3a^{3}b^{5}b_{22}))+a(3a^{2}(2+a(-3+b))(-1+b)b^{4}a_{2}^{3}
\\
& -3a^{2}b^{3}a_{2}^{2}((-1+a^{2}(-2+b)+b-b^{2}+a(3-2b+b^{2}))b_{1}+b(-4+5b)
\\
& +a^{2}(3-4b+2b^{2})+a(-2+7b-8b^{2}+b^{3}))b_{2}) \\
& -ba_{2}(-3(-1+b)b^{2}(3a^{4}-7a^{3}(1+b)+2b^{2}(1+b)-ab(5+9b+2b^{2}) \\
& +a^{2}(3+14b+6b^{2}))a_{11}+3(-1+b)b^{2}(6b^{3}+3a^{4}(2+b)-ab^{2}(17+8b)
\\
& +a^{2}b(15+26b)-a^{3}(4+24b+3b^{2}))a_{12}+a(3a(-1+b)b^{3}(a^{2}(-4+b)-3b
\\
& +a(3+4b-b^{2}))a_{22}-3a(-1+b)(a^{2}(-2+b)-3b+a(1+5b-2b^{2}))b_{1}^{2} \\
& +3a(a^{3}(1-b+b^{2})+b^{2}(3-5b+3b^{2})+a^{2}(1-5b+7b^{2}-4b^{3}) \\
& +ab(-2+b+b^{2}-b^{3}))b_{1}b_{2}+b(-3ab(a^{3}b+b(-2+3b)+a(-4+7b-4b^{2}) \\
& +a^{2}(5-8b+2b^{2}))b_{2}^{2}+(a-b)(-1+b)(3a(1+a(-2+b))b_{11} \\
& +3a(a+a^{2}(1-2b)-ab^{2}+b(-2+3b))b_{12}-b((a^{3}(11-7b)-3b^{2} \\
& +2ab(1+4b)+a^{2}(1-19b+7b^{2}))K+3a(a+a^{2}(-2+b)+2b-2ab)b_{22}))))) \\
&
+(-1+a)a(-3(a-b)^{2}(-1+b)b^{3}a_{111}-3(-2+3a-2b)(a-b)^{2}(-1+b)b^{3}a_{112}
\\
&
-6a^{3}b^{3}a_{122}+21a^{2}b^{4}a_{122}-24ab^{5}a_{122}-9a^{2}b^{5}a_{122}+6a^{3}b^{5}a_{122}+9b^{6}a_{122}
\\
&
+18ab^{6}a_{122}-12a^{2}b^{6}a_{122}-9b^{7}a_{122}+6ab^{7}a_{122}+3a^{3}b^{4}a_{222}-6a^{2}b^{5}a_{222}
\\
&
-3a^{3}b^{5}a_{222}+3ab^{6}a_{222}+6a^{2}b^{6}a_{222}-3ab^{7}a_{222}+7a^{3}b^{2}Kb_{1}-11a^{4}b^{2}Kb_{1}
\\
&
-14a^{2}b^{3}Kb_{1}+30a^{3}b^{3}Kb_{1}-a^{4}b^{3}Kb_{1}+7ab^{4}Kb_{1}-27a^{2}b^{4}Kb_{1}-a^{3}b^{4}Kb_{1}
\\
&
+8ab^{5}Kb_{1}+5a^{2}b^{5}Kb_{1}-3ab^{6}Kb_{1}-3ab^{2}a_{11}b_{1}+6a^{2}b^{2}a_{11}b_{1}+3b^{3}a_{11}b_{1}
\\
&
-3a^{2}b^{3}a_{11}b_{1}-6b^{4}a_{11}b_{1}-3ab^{4}a_{11}b_{1}+6b^{5}a_{11}b_{1}-6a^{2}b^{2}a_{12}b_{1}+3a^{3}b^{2}a_{12}b_{1}
\\
&
+3ab^{3}a_{12}b_{1}-3a^{2}b^{3}a_{12}b_{1}+3a^{3}b^{3}a_{12}b_{1}+3b^{4}a_{12}b_{1}+9ab^{4}a_{12}b_{1}-9a^{2}b^{4}a_{12}b_{1}
\\
&
-9b^{5}a_{12}b_{1}+6ab^{5}a_{12}b_{1}+3a^{2}b^{3}a_{22}b_{1}-6a^{3}b^{3}a_{22}b_{1}-3ab^{4}a_{22}b_{1}+6a^{2}b^{4}a_{22}b_{1}
\\
&
+3a^{3}b^{4}a_{22}b_{1}-3a^{2}b^{5}a_{22}b_{1}-3a^{2}b_{1}^{3}-3a^{3}b_{1}^{3}+6abb_{1}^{3}+12a^{2}bb_{1}^{3}+6a^{3}bb_{1}^{3}
\\
&
-18ab^{2}b_{1}^{3}-18a^{2}b^{2}b_{1}^{3}+18ab^{3}b_{1}^{3}-3a^{4}b^{2}Kb_{2}+5a^{3}b^{3}Kb_{2}+8a^{4}b^{3}Kb_{2}
\\
&
-a^{2}b^{4}Kb_{2}-27a^{3}b^{4}Kb_{2}+7a^{4}b^{4}Kb_{2}-ab^{5}Kb_{2}+30a^{2}b^{5}Kb_{2}-14a^{3}b^{5}Kb_{2}
\\
&
-11ab^{6}Kb_{2}+7a^{2}b^{6}Kb_{2}-3ab^{3}a_{11}b_{2}+3b^{4}a_{11}b_{2}+6ab^{4}a_{11}b_{2}-3a^{2}b^{4}a_{11}b_{2}
\\
&
-6b^{5}a_{11}b_{2}+3ab^{5}a_{11}b_{2}+6a^{2}b^{3}a_{12}b_{2}-9a^{3}b^{3}a_{12}b_{2}-9ab^{4}a_{12}b_{2}+9a^{2}b^{4}a_{12}b_{2}
\\
&
+3a^{3}b^{4}a_{12}b_{2}+3b^{5}a_{12}b_{2}-3ab^{5}a_{12}b_{2}+3a^{2}b^{5}a_{12}b_{2}+3b^{6}a_{12}b_{2}-6ab^{6}a_{12}b_{2}%
\end{array}%
\end{equation*}

\newpage

\begin{equation*}
\begin{array}{ll}
&
+6a^{3}b^{3}a_{22}b_{2}-3a^{2}b^{4}a_{22}b_{2}-6a^{3}b^{4}a_{22}b_{2}-3ab^{5}a_{22}b_{2}+3a^{3}b^{5}a_{22}b_{2}+6ab^{6}a_{22}b_{2}
\\
&
-3a^{2}b^{6}a_{22}b_{2}+15a^{3}b_{1}^{2}b_{2}+6a^{4}b_{1}^{2}b_{2}-39a^{2}bb_{1}^{2}b_{2}-36a^{3}bb_{1}^{2}b_{2}-12a^{4}bb_{1}^{2}b_{2}
\\
&
+30ab^{2}b_{1}^{2}b_{2}+69a^{2}b^{2}b_{1}^{2}b_{2}+36a^{3}b^{2}b_{1}^{2}b_{2}-48ab^{3}b_{1}^{2}b_{2}-39a^{2}b^{3}b_{1}^{2}b_{2}+18ab^{4}b_{1}^{2}b_{2}
\\
&
-18a^{4}b_{1}b_{2}^{2}+51a^{3}bb_{1}b_{2}^{2}+24a^{4}bb_{1}b_{2}^{2}-48a^{2}b^{2}b_{1}b_{2}^{2}-69a^{3}b^{2}b_{1}b_{2}^{2}+12ab^{3}b_{1}b_{2}^{2}
\\
&
+72a^{2}b^{3}b_{1}b_{2}^{2}+3a^{3}b^{3}b_{1}b_{2}^{2}-18ab^{4}b_{1}b_{2}^{2}-9a^{2}b^{4}b_{1}b_{2}^{2}-6a^{2}b^{3}b_{2}^{3}+9a^{2}b^{4}b_{2}^{3}
\\
&
-3a^{3}b^{4}b_{2}^{3}+3a^{2}bb_{1}b_{11}+6a^{3}bb_{1}b_{11}-3ab^{2}b_{1}b_{11}-24a^{2}b^{2}b_{1}b_{11}-9a^{3}b^{2}b_{1}b_{11}
\\
&
+18ab^{3}b_{1}b_{11}+27a^{2}b^{3}b_{1}b_{11}-18ab^{4}b_{1}b_{11}-6a^{3}bb_{2}b_{11}-6a^{4}bb_{2}b_{11}+15a^{2}b^{2}b_{2}b_{11}
\\
&
+27a^{3}b^{2}b_{2}b_{11}+6a^{4}b^{2}b_{2}b_{11}-9ab^{3}b_{2}b_{11}-42a^{2}b^{3}b_{2}b_{11}-18a^{3}b^{3}b_{2}b_{11}
\\
&
+21ab^{4}b_{2}b_{11}+21a^{2}b^{4}b_{2}b_{11}-9ab^{5}b_{2}b_{11}-15a^{3}bb_{1}b_{12}-6a^{4}bb_{1}b_{12}+39a^{2}b^{2}b_{1}b_{12}
\\
&
+36a^{3}b^{2}b_{1}b_{12}+12a^{4}b^{2}b_{1}b_{12}-24ab^{3}b_{1}b_{12}-66a^{2}b^{3}b_{1}b_{12}-39a^{3}b^{3}b_{1}b_{12}
\\
&
+36ab^{4}b_{1}b_{12}+45a^{2}b^{4}b_{1}b_{12}-18ab^{5}b_{1}b_{12}+18a^{4}bb_{2}b_{12}-51a^{3}b^{2}b_{2}b_{12}
\\
&
-24a^{4}b^{2}b_{2}b_{12}+45a^{2}b^{3}b_{2}b_{12}+78a^{3}b^{3}b_{2}b_{12}-12ab^{4}b_{2}b_{12}-72a^{2}b^{4}b_{2}b_{12}
\\
&
-9a^{3}b^{4}b_{2}b_{12}+18ab^{5}b_{2}b_{12}+9a^{2}b^{5}b_{2}b_{12}+9a^{4}bb_{1}b_{22}-27a^{3}b^{2}b_{1}b_{22}
\\
&
-9a^{4}b^{2}b_{1}b_{22}+24a^{2}b^{3}b_{1}b_{22}+30a^{3}b^{3}b_{1}b_{22}+3a^{4}b^{3}b_{1}b_{22}-6ab^{4}b_{1}b_{22}
\\
&
-27a^{2}b^{4}b_{1}b_{22}-9a^{3}b^{4}b_{1}b_{22}+6ab^{5}b_{1}b_{22}+6a^{2}b^{5}b_{1}b_{22}-9a^{3}b^{3}b_{2}b_{22}
\\
&
+9a^{2}b^{4}b_{2}b_{22}+12a^{3}b^{4}b_{2}b_{22}-3a^{4}b^{4}b_{2}b_{22}-12a^{2}b^{5}b_{2}b_{22}+3a^{3}b^{5}b_{2}b_{22}
\\
&
-3a^{3}b^{2}b_{111}+6a^{2}b^{3}b_{111}+3a^{3}b^{3}b_{111}-3ab^{4}b_{111}-6a^{2}b^{4}b_{111}+3ab^{5}b_{111}
\\
&
+6a^{3}b^{2}b_{112}+6a^{4}b^{2}b_{112}-12a^{2}b^{3}b_{112}-27a^{3}b^{3}b_{112}-6a^{4}b^{3}b_{112}+6ab^{4}b_{112}
\\
&
+36a^{2}b^{4}b_{112}+21a^{3}b^{4}b_{112}-15ab^{5}b_{112}-24a^{2}b^{5}b_{112}+9ab^{6}b_{112}-9a^{4}b^{2}b_{122}
\\
&
+24a^{3}b^{3}b_{122}+15a^{4}b^{3}b_{122}-21a^{2}b^{4}b_{122}-36a^{3}b^{4}b_{122}-6a^{4}b^{4}b_{122}+6ab^{5}b_{122}
\\
&
+27a^{2}b^{5}b_{122}+12a^{3}b^{5}b_{122}-6ab^{6}b_{122}-6a^{2}b^{6}b_{122}-3a^{4}b^{3}b_{222}+6a^{3}b^{4}b_{222}
\\
&
+3a^{4}b^{4}b_{222}-3a^{2}b^{5}b_{222}-6a^{3}b^{5}b_{222}+3a^{2}b^{6}b_{222}-15a^{4}b^{3}K_{1}+45a^{3}b^{4}K_{1}
\\
&
+15a^{4}b^{4}K_{1}-45a^{2}b^{5}K_{1}-45a^{3}b^{5}K_{1}+15ab^{6}K_{1}+45a^{2}b^{6}K_{1}-15ab^{7}K_{1}
\\
&
+15a^{4}b^{3}K_{2}-45a^{3}b^{4}K_{2}-15a^{4}b^{4}K_{2}+45a^{2}b^{5}K_{2}+45a^{3}b^{5}K_{2}-15ab^{6}K_{2}
\\
& -45a^{2}b^{6}K_{2}+15ab^{7}K_{2}));%
\end{array}%
\end{equation*}

\begin{equation*}
\begin{array}{ll}
j_{3}= & -(-1+b)b^{2}(6a^{3}-6a^{2}(1+b)-b(1+b)+2a(1+b)^{2})Ka_{1}^{2} \\
& +2a^{2}(-2-5b+5b^{2}+2b^{3}))Ka_{2}-ba_{1}(b(3a^{4}(-1+b)+3(-1+b)b^{2} \\
& -8a^{3}(-1+b^{2})-6ab(-1+b^{2})+a((b+a^{2}(7-12b)b-2b^{3} \\
& +a^{3}(-2+4b)+a(1-8b+7b^{2}+4b^{3}))Kb_{1}+b((-(-2+b)b^{2} \\
& +a^{3}(-6+4b)+a^{2}(3+7b-5b^{2})+a(1-6b-b^{2}+2b^{3}))Kb_{2} \\
& +(-1+b)(3a^{3}-b^{2}+ab(3+2b)-a^{2}(2+5b))(-K_{1}+bK_{2})))) \\
& -a(-a(-1+b)b^{2}(-6ab+a^{2}(1+b)+2b(1+b))Ka_{2}^{2} \\
& +aba_{2}((b+3b^{2}-6b^{3}+a^{3}(-1+2b)+ab(-6+7b+4b^{2}) \\
& -a^{2}(-2+b+5b^{2}))Kb_{1}+b((a^{3}(-3+2b)+b(-2+7b-3b^{2}) \\
& -a^{2}(-7+b+2b^{2})+a(-2-4b-b^{2}+2b^{3}))Kb_{2} \\
& +(-1+b)(a^{3}-3b^{2}+ab(5+2b)-a^{2}(2+3b))(-K_{1}+bK_{2}))) \\
& +(-1+a)(-(-1+b)b^{2}(-a+3a^{2}+b-4ab+b^{2})Ka_{11} \\
& -(-1+b)b^{2}(3a^{3}-3b^{2}+ab(7+4b)-a^{2}(4+7b))Ka_{12} \\
& +a((-1+b)b^{2}(3b^{2}+a^{2}(1+b)-ab(4+b))Ka_{22} \\
& +(a^{2}(-1+2b)+a(-1+4b-6b^{2})+2b(1-3b+3b^{2}))Kb_{1}^{2} \\
& +b_{1}((b^{2}(4-8b+3b^{2})-2ab(3-7b+4b^{2})+a^{2}(3-6b+4b^{2}))Kb_{2} \\
& +b(a(3-5b)b+a^{2}(-1+2b)+b^{2}(-2+3b))(-K_{1}+aK_{2}))%
\end{array}%
\end{equation*}

\begin{equation*}
\begin{array}{ll}
& +b(-b(2a^{2}+b^{2}+a(2-6b+b^{2}))Kb_{2}^{2} \\
& +b(a(5-3b)b+(-2+b)b^{2}+a^{2}(-3+2b))b_{2}(-K_{1}+aK_{2}) \\
& +(a-b)(-1+b)(-(1+a-3b)Kb_{11}+(a(3-4b)+b(-4+3b))Kb_{12} \\
& -b((a(-3+b)+b)Kb_{22}-(a-b)(-10K^{2}+5aK^{2}+5bK^{2}+K_{11} \\
& -(a+b)K_{12}+abK_{22})))))));%
\end{array}%
\end{equation*}

\begin{equation*}
\begin{array}{ll}
j_{4}= & -6(-1+b)b^{2}(6a^{3}-6a^{2}(1+b)-b(1+b)+2a(1+b)^{2})a_{1}^{3}b_{2}
\\
& -6ba_{1}^{2}(-(-1+b)b(-6a^{4}-b^{2}(5+2b)+a^{3}(16+13b)+ab(13+12b+4b^{2})
\\
& -a^{2}(10+23b+12b^{2}))a_{2}b_{2}+a((b+a^{2}(7-12b)b-2b^{3}+a^{3}(-2+4b)
\\
& +a(1-8b+7b^{2}+4b^{3}))b_{1}b_{2}+b(-(-1+a)(2a^{2}b+b(-1+2b) \\
& -a(-3+5b+b^{2}))b_{2}^{2}+(-1+b)(6a^{3}-6a^{2}(1+b)-b(1+b) \\
& +2a(1+b)^{2})(bK-b_{12}))))+2a_{1}(-3(-1+b)b^{2}(-6b^{3}-3a^{4}(2+b) \\
& +ab^{2}(17+8b)-a^{2}b(16+23b)+a^{3}(4+24b+b^{2}))a_{2}^{2}b_{2} \\
& -3aba_{2}((-1+a)(2(2-3b)b^{2}+a^{3}(-2+4b)+a^{2}(1+5b-12b^{2}) \\
& +ab(-8+9b+5b^{2}))b_{1}b_{2}+b((a^{4}(-3+b)-b^{2}(1+b)+ab(2+5b+b^{2}) \\
& +a^{3}(5-b+4b^{2})-a^{2}(3+b+7b^{2}+b^{3}))b_{2}^{2}+(-1+b)(3a^{4}+3b^{2}
\\
& -8a^{3}(1+b)-6ab(1+b)+2a^{2}(2+7b+2b^{2}))(bK-b_{12}))) \\
& +a(-3a(2a^{2}(1-3b)b+a^{3}(-1+2b)-2b(1-3b+3b^{2}) \\
& +a(1-2b+6b^{3}))b_{1}^{2}b_{2}+ab_{1}(3(a^{3}(-3+4b)+b^{2}(5-8b+b^{2}) \\
& +a^{2}(3-3b^{2}-4b^{3})+ab(-5+2b+7b^{2}+b^{3}))b_{2}^{2} \\
& +b(b(3a^{3}(-5+3b)+a^{2}(14+b^{2})-b(3-7b+b^{2}) \\
& +a(-3+3b-14b^{2}+2b^{3}))K+3(b+a^{2}(7-12b)b-2b^{3} \\
& +a^{3}(-2+4b)+a(1-8b+7b^{2}+4b^{3}))b_{12}-3(-1+b)(3a^{3}-b^{2}+ab(3+2b)
\\
& -a^{2}(2+5b))b_{22}))+b(3a^{2}(1+a^{2}+a(-1+b)-2b)(-1+b)bb_{2}^{3} \\
&
+b_{2}(9a^{3}bK+5a^{4}bK-3a^{5}bK-20a^{2}b^{2}K-34a^{3}b^{2}K+9a^{4}b^{2}K+3a^{5}b^{2}K
\\
&
+8ab^{3}K+54a^{2}b^{3}K+9a^{3}b^{3}K-8a^{4}b^{3}K-22ab^{4}K-23a^{2}b^{4}K+a^{3}b^{4}K
\\
& +11ab^{5}K+a^{2}b^{5}K+3(-1+b)b(6a^{3}-b(1+2b)-3a^{2}(2+3b) \\
& +a(1+8b+3b^{2}))a_{11}+3(-1+b)b(6a^{4}+b^{2}(5+2b)-3a^{3}(4+5b) \\
& -ab(13+12b+4b^{2})+a^{2}(8+22b+13b^{2}))a_{12}-6a^{3}ba_{22}+6a^{4}ba_{22}
\\
&
+24a^{2}b^{2}a_{22}-21a^{3}b^{2}a_{22}-27ab^{3}a_{22}+6a^{2}b^{3}a_{22}+18a^{3}b^{3}a_{22}-6a^{4}b^{3}a_{22}
\\
&
+9b^{4}a_{22}+18ab^{4}a_{22}-27a^{2}b^{4}a_{22}+9a^{3}b^{4}a_{22}-9b^{5}a_{22}+9ab^{5}a_{22}
\\
&
-3a^{2}b^{5}a_{22}+3a^{2}b_{11}-3a^{4}b_{11}-3abb_{11}-12a^{2}bb_{11}+12a^{3}bb_{11}+3a^{4}bb_{11}
\\
&
+12ab^{2}b_{11}-12a^{3}b^{2}b_{11}-9ab^{3}b_{11}+9a^{2}b^{3}b_{11}-9a^{3}b_{12}+9a^{4}b_{12}+24a^{2}bb_{12}
\\
&
-9a^{3}bb_{12}-12a^{4}bb_{12}-12ab^{2}b_{12}-21a^{2}b^{2}b_{12}+24a^{3}b^{2}b_{12}-3a^{4}b^{2}b_{12}
\\
&
+18ab^{3}b_{12}-12a^{2}b^{3}b_{12}+9a^{3}b^{3}b_{12}-3ab^{4}b_{12}-3a^{2}b^{4}b_{12}-3a^{3}bb_{22}
\\
&
+3a^{2}b^{2}b_{22}+9a^{3}b^{2}b_{22}-3a^{4}b^{2}b_{22}-9a^{2}b^{3}b_{22}-3a^{3}b^{3}b_{22}+3a^{4}b^{3}b_{22}
\\
& +6a^{2}b^{4}b_{22}-3a^{3}b^{4}b_{22})+3a(-1+b)b(3a^{3}-b^{2}+ab(3+2b) \\
& -a^{2}(2+5b))(-b_{112}+b_{122}+b(K_{1}-K_{2}))))) \\
&
-a(-6a^{2}(2+a(-3+b))(-1+b)b^{3}a_{2}^{3}b_{2}-6a^{2}ba_{2}^{2}((a^{2}(-1+2b)
\\
& -ab(-3+4b+b^{2})+b(-1-b+3b^{2}))b_{1}b_{2}+b(-(b(-4+5b) \\
& +a^{2}(3-4b+2b^{2})+a(-2+7b-8b^{2}+b^{3}))b_{2}^{2}+(-1+b)(-6ab+a^{2}(1+b)
\\
& +2b(1+b))(bK-b_{12})))+2a_{2}(3a^{2}(2a^{2}(1-3b)b+a^{3}(-1+2b) \\
&
-2b(1-3b+3b^{2})+a(1-2b+6b^{3}))b_{1}^{2}b_{2}-a^{2}b_{1}(-3(a^{3}(3-5b+2b^{2})%
\end{array}%
\end{equation*}

\newpage

\begin{equation*}
\begin{array}{ll}
&
+b^{2}(-5+5b+3b^{2})-ab(-6+4b+7b^{2}+b^{3})-a^{2}(3+2b-11b^{2}+3b^{3}))b_{2}^{2}
\\
& +b(b(a^{3}(-4+b)+a^{2}(8+10b-6b^{2})-3b(1-4b+b^{2})+ab(-17+2b^{2}))K \\
& +3(b+3b^{2}-6b^{3}+a^{3}(-1+2b)+ab(-6+7b+4b^{2})-a^{2}(-2+b+5b^{2}))b_{12}
\\
& -3(-1+b)(a^{3}-3b^{2}+ab(5+2b)-a^{2}(2+3b))b_{22})) \\
& +b(-3a^{2}b(a^{3}b+b(-2+3b)+a(-4+7b-4b^{2})+a^{2}(5-8b+2b^{2}))b_{2}^{3}
\\
& +b_{2}(-3(-1+b)b(3a^{4}-7a^{3}(1+b)+2b^{2}(1+b)-ab(5+9b+2b^{2}) \\
& +a^{2}(3+14b+6b^{2}))a_{11}+3(-1+b)b(6b^{3}+3a^{4}(2+b)-ab^{2}(17+8b) \\
&
+a^{2}b(15+26b)-a^{3}(4+24b+3b^{2}))a_{12}+a(-a^{3}bK+12a^{4}bK-17a^{3}b^{2}K
\\
&
-25a^{4}b^{2}K-14ab^{3}K+26a^{2}b^{3}K+41a^{3}b^{3}K+10a^{4}b^{3}K+3b^{4}K+12ab^{4}K
\\
&
-48a^{2}b^{4}K-11a^{3}b^{4}K-3b^{5}K+8ab^{5}K+7a^{2}b^{5}K+3a(-1+b)b^{2}(a^{2}(-4+b)
\\
& -3b+a(3+4b-b^{2}))a_{22}-3(-1+a)a(-1+b)(a+a^{2}-4ab+b(-1+3b))b_{11} \\
&
+9a^{3}b_{12}-9a^{4}b_{12}-15a^{2}bb_{12}-3a^{3}bb_{12}+21a^{4}bb_{12}+18ab^{2}b_{12}-3a^{2}b^{2}b_{12}
\\
&
-18a^{3}b^{2}b_{12}-9a^{4}b^{2}b_{12}-15ab^{3}b_{12}+30a^{2}b^{3}b_{12}-9ab^{4}b_{12}+3a^{2}b^{4}b_{12}+3a^{3}bb_{22}
\\
&
-6a^{4}bb_{22}+3a^{2}b^{2}b_{22}-3a^{3}b^{2}b_{22}+9a^{4}b^{2}b_{22}-6ab^{3}b_{22}+3a^{2}b^{3}b_{22}-3a^{3}b^{3}b_{22}
\\
&
-3a^{4}b^{3}b_{22}+6ab^{4}b_{22}-6a^{2}b^{4}b_{22}+3a^{3}b^{4}b_{22}))-3a^{2}(-1+b)b(a^{3}-3b^{2}
\\
& +ab(5+2b)-a^{2}(2+3b))(-b_{112}+b_{122}+b(K_{1}-K_{2})))) \\
&
+(-1+a)a(30a^{3}b^{3}K^{2}-30a^{4}b^{3}K^{2}-60a^{2}b^{4}K^{2}+30a^{3}b^{4}K^{2}
\\
&
+30a^{4}b^{4}K^{2}+30ab^{5}K^{2}+30a^{2}b^{5}K^{2}-60a^{3}b^{5}K^{2}-30ab^{6}K^{2}+30a^{2}b^{6}K^{2}
\\
&
-6a^{3}b^{3}Ka_{22}+24a^{2}b^{4}Ka_{22}-18ab^{5}Ka_{22}-18a^{2}b^{5}Ka_{22}+6a^{3}b^{5}Ka_{22}
\\
&
+18ab^{6}Ka_{22}-6a^{2}b^{6}Ka_{22}-6a^{2}bKb_{1}^{2}+8a^{3}bKb_{1}^{2}+12ab^{2}Kb_{1}^{2}-18a^{2}b^{2}Kb_{1}^{2}
\\
&
-16a^{3}b^{2}Kb_{1}^{2}-8ab^{3}Kb_{1}^{2}+34a^{2}b^{3}Kb_{1}^{2}-6ab^{4}Kb_{1}^{2}-6a^{2}b^{2}a_{111}b_{2}+12ab^{3}a_{111}b_{2}
\\
&
+6a^{2}b^{3}a_{111}b_{2}-6b^{4}a_{111}b_{2}-12ab^{4}a_{111}b_{2}+6b^{5}a_{111}b_{2}+12a^{2}b^{2}a_{112}b_{2}
\\
&
-18a^{3}b^{2}a_{112}b_{2}-24ab^{3}a_{112}b_{2}+36a^{2}b^{3}a_{112}b_{2}+18a^{3}b^{3}a_{112}b_{2}+12b^{4}a_{112}b_{2}
\\
&
-18ab^{4}a_{112}b_{2}-48a^{2}b^{4}a_{112}b_{2}+42ab^{5}a_{112}b_{2}-12b^{6}a_{112}b_{2}+12a^{3}b^{2}a_{122}b_{2}
\\
&
-42a^{2}b^{3}a_{122}b_{2}+48ab^{4}a_{122}b_{2}+18a^{2}b^{4}a_{122}b_{2}-12a^{3}b^{4}a_{122}b_{2}-18b^{5}a_{122}b_{2}
\\
&
-36ab^{5}a_{122}b_{2}+24a^{2}b^{5}a_{122}b_{2}+18b^{6}a_{122}b_{2}-12ab^{6}a_{122}b_{2}-6a^{3}b^{3}a_{222}b_{2}
\\
&
+12a^{2}b^{4}a_{222}b_{2}+6a^{3}b^{4}a_{222}b_{2}-6ab^{5}a_{222}b_{2}-12a^{2}b^{5}a_{222}b_{2}+6ab^{6}a_{222}b_{2}
\\
&
+20a^{3}bKb_{1}b_{2}-22a^{4}bKb_{1}b_{2}-42a^{2}b^{2}Kb_{1}b_{2}+10a^{3}b^{2}Kb_{1}b_{2}+44a^{4}b^{2}Kb_{1}b_{2}
\\
&
+28ab^{3}Kb_{1}b_{2}+68a^{2}b^{3}Kb_{1}b_{2}-84a^{3}b^{3}Kb_{1}b_{2}-56ab^{4}Kb_{1}b_{2}+28a^{2}b^{4}Kb_{1}b_{2}
\\
&
+6ab^{5}Kb_{1}b_{2}-6a^{3}ba_{22}b_{1}b_{2}+18a^{2}b^{2}a_{22}b_{1}b_{2}+12a^{3}b^{2}a_{22}b_{1}b_{2}-12ab^{3}a_{22}b_{1}b_{2}
\\
&
-30a^{2}b^{3}a_{22}b_{1}b_{2}+18ab^{4}a_{22}b_{1}b_{2}+4a^{3}b^{2}Kb_{2}^{2}-8a^{4}b^{2}Kb_{2}^{2}-18a^{2}b^{3}Kb_{2}^{2}
\\
&
+36a^{3}b^{3}Kb_{2}^{2}-14a^{4}b^{3}Kb_{2}^{2}+2ab^{4}Kb_{2}^{2}-32a^{2}b^{4}Kb_{2}^{2}+22a^{3}b^{4}Kb_{2}^{2}+22ab^{5}Kb_{2}^{2}
\\
&
-14a^{2}b^{5}Kb_{2}^{2}-12a^{3}b^{2}a_{22}b_{2}^{2}+6a^{2}b^{3}a_{22}b_{2}^{2}+12a^{3}b^{3}a_{22}b_{2}^{2}+6ab^{4}a_{22}b_{2}^{2}
\\
&
-6a^{3}b^{4}a_{22}b_{2}^{2}-12ab^{5}a_{22}b_{2}^{2}+6a^{2}b^{5}a_{22}b_{2}^{2}+6a^{2}b_{1}^{2}b_{2}^{2}+6a^{3}b_{1}^{2}b_{2}^{2}-12abb_{1}^{2}b_{2}^{2}
\\
&
-24a^{2}bb_{1}^{2}b_{2}^{2}-12a^{3}bb_{1}^{2}b_{2}^{2}+36ab^{2}b_{1}^{2}b_{2}^{2}+36a^{2}b^{2}b_{1}^{2}b_{2}^{2}-36ab^{3}b_{1}^{2}b_{2}^{2}-18a^{3}b_{1}b_{2}^{3}
\\
&
+36a^{2}bb_{1}b_{2}^{3}+30a^{3}bb_{1}b_{2}^{3}-24ab^{2}b_{1}b_{2}^{3}-66a^{2}b^{2}b_{1}b_{2}^{3}-12a^{3}b^{2}b_{1}b_{2}^{3}+36ab^{3}b_{1}b_{2}^{3}
\\
&
+18a^{2}b^{3}b_{1}b_{2}^{3}+12a^{2}b^{2}b_{2}^{4}-18a^{2}b^{3}b_{2}^{4}+6a^{3}b^{3}b_{2}^{4}+6a^{2}b^{2}Kb_{11}-18a^{3}b^{2}Kb_{11}
\\
&
-6ab^{3}Kb_{11}+18a^{2}b^{3}Kb_{11}+18a^{3}b^{3}Kb_{11}-24a^{2}b^{4}Kb_{11}+6ab^{5}Kb_{11}-6a^{2}bb_{2}^{2}b_{11}
\\
&
-6a^{3}bb_{2}^{2}b_{11}+6ab^{2}b_{2}^{2}b_{11}+30a^{2}b^{2}b_{2}^{2}b_{11}+6a^{3}b^{2}b_{2}^{2}b_{11}-24ab^{3}b_{2}^{2}b_{11}-24a^{2}b^{3}b_{2}^{2}b_{11}
\\
& +18ab^{4}b_{2}^{2}b_{11}-6(a-b)ba_{11}(-(a-2ab+b(-2+3b))b_{1}b_{2} \\
& -b((1+(-2+a)b)b_{2}^{2}+(-1+b)(1-3a+b)(bK-b_{12}))) \\
& -6(a-b)ba_{12}(2(a-3ab^{2}+a^{2}(-1+2b)+b(-2+3b))b_{1}b_{2} \\
& +b((a^{2}(-3+b)-b(1+b)+2a(1+b^{2}))b_{2}^{2} \\
& +(-1+b)(3a^{2}+3b-4a(1+b))(bK-b_{12})))%
\end{array}%
\end{equation*}

\newpage

\begin{equation*}
\begin{array}{ll}
&
-36a^{3}b^{2}Kb_{12}+42a^{4}b^{2}Kb_{12}+78a^{2}b^{3}Kb_{12}-42a^{3}b^{3}Kb_{12}-42a^{4}b^{3}Kb_{12}
\\
&
-42ab^{4}Kb_{12}-42a^{2}b^{4}Kb_{12}+78a^{3}b^{4}Kb_{12}+42ab^{5}Kb_{12}-36a^{2}b^{5}Kb_{12}
\\
&
+6a^{3}b^{2}a_{22}b_{12}-24a^{2}b^{3}a_{22}b_{12}+18ab^{4}a_{22}b_{12}+18a^{2}b^{4}a_{22}b_{12}-6a^{3}b^{4}a_{22}b_{12}
\\
&
-18ab^{5}a_{22}b_{12}+6a^{2}b^{5}a_{22}b_{12}+6a^{2}b_{1}^{2}b_{12}+6a^{3}b_{1}^{2}b_{12}-12abb_{1}^{2}b_{12}
\\
&
-24a^{2}bb_{1}^{2}b_{12}-12a^{3}bb_{1}^{2}b_{12}+36ab^{2}b_{1}^{2}b_{12}+36a^{2}b^{2}b_{1}^{2}b_{12}-36ab^{3}b_{1}^{2}b_{12}
\\
&
-18a^{3}b_{1}b_{2}b_{12}+36a^{2}bb_{1}b_{2}b_{12}+18a^{3}bb_{1}b_{2}b_{12}-24ab^{2}b_{1}b_{2}b_{12}-30a^{2}b^{2}b_{1}b_{2}b_{12}
\\
&
+12a^{3}b^{2}b_{1}b_{2}b_{12}+12ab^{3}b_{1}b_{2}b_{12}-42a^{2}b^{3}b_{1}b_{2}b_{12}+36ab^{4}b_{1}b_{2}b_{12}+18a^{3}bb_{2}^{2}b_{12}
\\
&
-30a^{2}b^{2}b_{2}^{2}b_{12}-48a^{3}b^{2}b_{2}^{2}b_{12}+24ab^{3}b_{2}^{2}b_{12}+66a^{2}b^{3}b_{2}^{2}b_{12}+24a^{3}b^{3}b_{2}^{2}b_{12}
\\
&
-36ab^{4}b_{2}^{2}b_{12}-18a^{2}b^{4}b_{2}^{2}b_{12}-6a^{2}bb_{11}b_{12}-6a^{3}bb_{11}b_{12}+6ab^{2}b_{11}b_{12}
\\
&
+30a^{2}b^{2}b_{11}b_{12}+6a^{3}b^{2}b_{11}b_{12}-24ab^{3}b_{11}b_{12}-24a^{2}b^{3}b_{11}b_{12}+18ab^{4}b_{11}b_{12}
\\
&
+18a^{3}bb_{12}^{2}-42a^{2}b^{2}b_{12}^{2}-24a^{3}b^{2}b_{12}^{2}+24ab^{3}b_{12}^{2}+48a^{2}b^{3}b_{12}^{2}+6a^{3}b^{3}b_{12}^{2}
\\
&
-24ab^{4}b_{12}^{2}-6a^{2}b^{4}b_{12}^{2}+6a^{4}b^{2}Kb_{22}-24a^{3}b^{3}Kb_{22}+18a^{2}b^{4}Kb_{22}+18a^{3}b^{4}Kb_{22}
\\
&
-6a^{4}b^{4}Kb_{22}-18a^{2}b^{5}Kb_{22}+6a^{3}b^{5}Kb_{22}-6a^{3}b_{1}^{2}b_{22}-6a^{4}b_{1}^{2}b_{22}+18a^{2}bb_{1}^{2}b_{22}
\\
&
+18a^{3}bb_{1}^{2}b_{22}+12a^{4}bb_{1}^{2}b_{22}-12ab^{2}b_{1}^{2}b_{22}-24a^{2}b^{2}b_{1}^{2}b_{22}-24a^{3}b^{2}b_{1}^{2}b_{22}
\\
&
+12ab^{3}b_{1}^{2}b_{22}+12a^{2}b^{3}b_{1}^{2}b_{22}+18a^{4}b_{1}b_{2}b_{22}-36a^{3}bb_{1}b_{2}b_{22}-18a^{4}bb_{1}b_{2}b_{22}
\\
&
+6a^{2}b^{2}b_{1}b_{2}b_{22}+36a^{3}b^{2}b_{1}b_{2}b_{22}-12a^{4}b^{2}b_{1}b_{2}b_{22}+12ab^{3}b_{1}b_{2}b_{22}-6a^{2}b^{3}b_{1}b_{2}b_{22}
\\
&
+24a^{3}b^{3}b_{1}b_{2}b_{22}-12ab^{4}b_{1}b_{2}b_{22}-12a^{2}b^{4}b_{1}b_{2}b_{22}+18a^{3}b^{2}b_{2}^{2}b_{22}-18a^{2}b^{3}b_{2}^{2}b_{22}
\\
&
-24a^{3}b^{3}b_{2}^{2}b_{22}+6a^{4}b^{3}b_{2}^{2}b_{22}+24a^{2}b^{4}b_{2}^{2}b_{22}-6a^{3}b^{4}b_{2}^{2}b_{22}+6a^{3}bb_{11}b_{22}+6a^{4}bb_{11}b_{22}
\\
&
-12a^{2}b^{2}b_{11}b_{22}-18a^{3}b^{2}b_{11}b_{22}-6a^{4}b^{2}b_{11}b_{22}+6ab^{3}b_{11}b_{22}+18a^{2}b^{3}b_{11}b_{22}
\\
&
+12a^{3}b^{3}b_{11}b_{22}-6ab^{4}b_{11}b_{22}-6a^{2}b^{4}b_{11}b_{22}-18a^{4}bb_{12}b_{22}+48a^{3}b^{2}b_{12}b_{22}
\\
&
+12a^{4}b^{2}b_{12}b_{22}-30a^{2}b^{3}b_{12}b_{22}-42a^{3}b^{3}b_{12}b_{22}+6a^{4}b^{3}b_{12}b_{22}+30a^{2}b^{4}b_{12}b_{22}
\\
&
-6a^{3}b^{4}b_{12}b_{22}-6a^{3}bb_{1}b_{112}+18a^{2}b^{2}b_{1}b_{112}+12a^{3}b^{2}b_{1}b_{112}-12ab^{3}b_{1}b_{112}
\\
&
-30a^{2}b^{3}b_{1}b_{112}+18ab^{4}b_{1}b_{112}+6a^{3}b^{2}b_{2}b_{112}-18a^{2}b^{3}b_{2}b_{112}-12a^{3}b^{3}b_{2}b_{112}
\\
&
+12ab^{4}b_{2}b_{112}+30a^{2}b^{4}b_{2}b_{112}-18ab^{5}b_{2}b_{112}+6a^{3}bb_{1}b_{122}+6a^{4}bb_{1}b_{122}
\\
&
-18a^{2}b^{2}b_{1}b_{122}-18a^{3}b^{2}b_{1}b_{122}-12a^{4}b^{2}b_{1}b_{122}+12ab^{3}b_{1}b_{122}+24a^{2}b^{3}b_{1}b_{122}
\\
&
+24a^{3}b^{3}b_{1}b_{122}-12ab^{4}b_{1}b_{122}-12a^{2}b^{4}b_{1}b_{122}-6a^{3}b^{2}b_{2}b_{122}-6a^{4}b^{2}b_{2}b_{122}
\\
&
+18a^{2}b^{3}b_{2}b_{122}+18a^{3}b^{3}b_{2}b_{122}+12a^{4}b^{3}b_{2}b_{122}-12ab^{4}b_{2}b_{122}-24a^{2}b^{4}b_{2}b_{122}
\\
&
-24a^{3}b^{4}b_{2}b_{122}+12ab^{5}b_{2}b_{122}+12a^{2}b^{5}b_{2}b_{122}-6a^{4}bb_{1}b_{222}+12a^{3}b^{2}b_{1}b_{222}
\\
&
+6a^{4}b^{2}b_{1}b_{222}-6a^{2}b^{3}b_{1}b_{222}-12a^{3}b^{3}b_{1}b_{222}+6a^{2}b^{4}b_{1}b_{222}+6a^{4}b^{2}b_{2}b_{222}
\\
&
-12a^{3}b^{3}b_{2}b_{222}-6a^{4}b^{3}b_{2}b_{222}+6a^{2}b^{4}b_{2}b_{222}+12a^{3}b^{4}b_{2}b_{222}-6a^{2}b^{5}b_{2}b_{222}
\\
&
+6a^{3}b^{2}b_{1112}-12a^{2}b^{3}b_{1112}-6a^{3}b^{3}b_{1112}+6ab^{4}b_{1112}+12a^{2}b^{4}b_{1112}
\\
&
-6ab^{5}b_{1112}-6a^{3}b^{2}b_{1122}-6a^{4}b^{2}b_{1122}+12a^{2}b^{3}b_{1122}+18a^{3}b^{3}b_{1122}
\\
&
+6a^{4}b^{3}b_{1122}-6ab^{4}b_{1122}-18a^{2}b^{4}b_{1122}-12a^{3}b^{4}b_{1122}+6ab^{5}b_{1122}+6a^{2}b^{5}b_{1122}
\\
&
+6a^{4}b^{2}b_{1222}-12a^{3}b^{3}b_{1222}-6a^{4}b^{3}b_{1222}+6a^{2}b^{4}b_{1222}+12a^{3}b^{4}b_{1222}
\\
&
-6a^{2}b^{5}b_{1222}-15a^{3}b^{2}b_{1}K_{1}+24a^{2}b^{3}b_{1}K_{1}+9a^{3}b^{3}b_{1}K_{1}-9ab^{4}b_{1}K_{1}-12a^{2}b^{4}b_{1}K_{1}
\\
&
+3ab^{5}b_{1}K_{1}-a^{3}b^{2}b_{2}K_{1}+29a^{4}b^{2}b_{2}K_{1}+2a^{2}b^{3}b_{2}K_{1}-75a^{3}b^{3}b_{2}K_{1}-29a^{4}b^{3}b_{2}K_{1}
\\
&
-ab^{4}b_{2}K_{1}+69a^{2}b^{4}b_{2}K_{1}+82a^{3}b^{4}b_{2}K_{1}-23ab^{5}b_{2}K_{1}-83a^{2}b^{5}b_{2}K_{1}+30ab^{6}b_{2}K_{1}
\\
&
+7a^{3}b^{2}b_{1}K_{2}+7a^{4}b^{2}b_{1}K_{2}-8a^{2}b^{3}b_{1}K_{2}-21a^{3}b^{3}b_{1}K_{2}-a^{4}b^{3}b_{1}K_{2}+ab^{4}b_{1}K_{2}
\\
&
+15a^{2}b^{4}b_{1}K_{2}+2a^{3}b^{4}b_{1}K_{2}-ab^{5}b_{1}K_{2}-a^{2}b^{5}b_{1}K_{2}-27a^{4}b^{2}b_{2}K_{2}+78a^{3}b^{3}b_{2}K_{2}
\\
&
+21a^{4}b^{3}b_{2}K_{2}-81a^{2}b^{4}b_{2}K_{2}-66a^{3}b^{4}b_{2}K_{2}+30ab^{5}b_{2}K_{2}+75a^{2}b^{5}b_{2}K_{2}
\\
&
-30ab^{6}b_{2}K_{2}-6a^{3}b^{3}K_{11}+12a^{2}b^{4}K_{11}+6a^{3}b^{4}K_{11}-6ab^{5}K_{11}-12a^{2}b^{5}K_{11}
\\
&
+6ab^{6}K_{11}+6a^{3}b^{3}K_{12}+6a^{4}b^{3}K_{12}-12a^{2}b^{4}K_{12}-18a^{3}b^{4}K_{12}-6a^{4}b^{4}K_{12}
\\
&
+6ab^{5}K_{12}+18a^{2}b^{5}K_{12}+12a^{3}b^{5}K_{12}-6ab^{6}K_{12}-6a^{2}b^{6}K_{12}-6a^{4}b^{3}K_{22}
\\
&
+12a^{3}b^{4}K_{22}+6a^{4}b^{4}K_{22}-6a^{2}b^{5}K_{22}-12a^{3}b^{5}K_{22}+6a^{2}b^{6}K_{22}));%
\end{array}%
\end{equation*}

\newpage

\begin{equation*}
\begin{array}{ll}
j_{5}= & -6(-1+b)b^{2}(6a^{3}-6a^{2}(1+b)-b(1+b)+2a(1+b)^{2})a_{1}^{3}a_{2}
\\
& -2ba_{1}^{2}(-3(-1+b)b(-b^{2}(5+2b)+a^{3}(10+7b)+ab(12+11b+4b^{2}) \\
& -a^{2}(8+19b+10b^{2}))a_{2}^{2}-a(-1+b)b(a(3a^{3}+a^{2}(4-17b) \\
& +(3-4b)b+a(-6+9b+8b^{2}))K+3(6a^{3}-6a^{2}(1+b)-b(1+b) \\
& +2a(1+b)^{2})a_{12}-3(2a^{3}(1+b)-b^{2}(1+b)-2a^{2}(1+2b+2b^{2}) \\
& +ab(3+3b+2b^{2}))a_{22})+3aa_{2}((b+a^{2}(7-12b)b-2b^{3}+a^{3}(-2+4b) \\
& +a(1-8b+7b^{2}+4b^{3}))b_{1}+b(-(-2+b)b^{2}+a^{3}(-6+4b) \\
& +a^{2}(3+7b-5b^{2})+a(1-6b-b^{2}+2b^{3}))b_{2})) \\
& -a_{1}(-6(-1+b)b^{3}(a^{3}(-13+b)+6b^{2}-2ab(7+4b)+2a^{2}(5+9b))a_{2}^{3}
\\
& +6aba_{2}^{2}((a^{3}(3-6b)+2b^{2}(-2+3b)+ab(7-5b-9b^{2}) \\
& +a^{2}(-2-5b+14b^{2}+b^{3}))b_{1}+b(2(-2+b)b^{2}+ab(8+3b-4b^{2}) \\
& +a^{3}(10-8b+b^{2})+a^{2}(-8-5b+4b^{2}+b^{3}))b_{2}) \\
&
+2aa_{2}(-21a^{3}b^{2}K+20a^{4}b^{2}K+35a^{2}b^{3}K+14a^{3}b^{3}K-39a^{4}b^{3}K
\\
&
-17ab^{4}K-70a^{2}b^{4}K+50a^{3}b^{4}K+19a^{4}b^{4}K+39ab^{5}K+14a^{2}b^{5}K
\\
& -43a^{3}b^{5}K-22ab^{6}K+21a^{2}b^{6}K-3(-1+b)b^{2}(6a^{3}-b(1+2b) \\
& -3a^{2}(2+3b)+a(1+8b+3b^{2}))a_{11}+3(-1+b)b^{2}(-2b^{2}(4+b) \\
&
+2a^{3}(7+4b)+ab(19+15b+4b^{2})-a^{2}(12+27b+11b^{2}))a_{12}-27a^{2}b^{3}a_{22}
\\
&
+30a^{3}b^{3}a_{22}+45ab^{4}a_{22}-21a^{2}b^{4}a_{22}-33a^{3}b^{4}a_{22}-18b^{5}a_{22}-27ab^{5}a_{22}
\\
&
+51a^{2}b^{5}a_{22}+3a^{3}b^{5}a_{22}+18b^{6}a_{22}-18ab^{6}a_{22}-3a^{2}b^{6}a_{22}+3a^{2}b_{1}^{2}-3a^{4}b_{1}^{2}
\\
&
-6abb_{1}^{2}-6a^{2}bb_{1}^{2}+6a^{3}bb_{1}^{2}+6a^{4}bb_{1}^{2}+18ab^{2}b_{1}^{2}-18a^{3}b^{2}b_{1}^{2}-18ab^{3}b_{1}^{2}
\\
&
+18a^{2}b^{3}b_{1}^{2}-9a^{3}b_{1}b_{2}+9a^{4}b_{1}b_{2}+18a^{2}bb_{1}b_{2}-18a^{4}bb_{1}b_{2}-12ab^{2}b_{1}b_{2}
\\
&
-30a^{2}b^{2}b_{1}b_{2}+30a^{3}b^{2}b_{1}b_{2}+12a^{4}b^{2}b_{1}b_{2}+24ab^{3}b_{1}b_{2}-24a^{3}b^{3}b_{1}b_{2}
\\
&
-9ab^{4}b_{1}b_{2}+9a^{2}b^{4}b_{1}b_{2}+6a^{2}b^{2}b_{2}^{2}-6a^{4}b^{2}b_{2}^{2}-18a^{2}b^{3}b_{2}^{2}+18a^{3}b^{3}b_{2}^{2}
\\
&
+3ab^{4}b_{2}^{2}-3a^{3}b^{4}b_{2}^{2}-3a^{2}bb_{11}+3a^{4}bb_{11}+3ab^{2}b_{11}+12a^{2}b^{2}b_{11}-12a^{3}b^{2}b_{11}
\\
&
-3a^{4}b^{2}b_{11}-12ab^{3}b_{11}+12a^{3}b^{3}b_{11}+9ab^{4}b_{11}-9a^{2}b^{4}b_{11}+9a^{3}bb_{12}-9a^{4}bb_{12}
\\
&
-21a^{2}b^{2}b_{12}+21a^{4}b^{2}b_{12}+12ab^{3}b_{12}+30a^{2}b^{3}b_{12}-30a^{3}b^{3}b_{12}-12a^{4}b^{3}b_{12}
\\
&
-21ab^{4}b_{12}+21a^{3}b^{4}b_{12}+9ab^{5}b_{12}-9a^{2}b^{5}b_{12}+9a^{3}b^{2}b_{22}-9a^{4}b^{2}b_{22}
\\
&
-12a^{2}b^{3}b_{22}+12a^{4}b^{3}b_{22}+3ab^{4}b_{22}+12a^{2}b^{4}b_{22}-12a^{3}b^{4}b_{22}-3a^{4}b^{4}b_{22}
\\
& -3ab^{5}b_{22}+3a^{3}b^{5}b_{22})+a^{2}b(6(-1+b)b(3a^{3}-b^{2}+ab(3+2b) \\
& -a^{2}(2+5b))a_{112}-6b(b^{2}-b^{4}+2a^{3}(-1+b^{2})+a^{2}(2+2b-4b^{3}) \\
&
+ab(-3+b^{2}+2b^{3}))a_{122}+6a^{2}b^{2}a_{222}-6a^{3}b^{2}a_{222}-12ab^{3}a_{222}+6a^{2}b^{3}a_{222}
\\
&
+6a^{3}b^{3}a_{222}+6b^{4}a_{222}+6ab^{4}a_{222}-12a^{2}b^{4}a_{222}-6b^{5}a_{222}+6ab^{5}a_{222}
\\
&
+6a^{2}Kb_{1}-14a^{3}Kb_{1}+2a^{4}Kb_{1}+6abKb_{1}-6a^{2}bKb_{1}+28a^{3}bKb_{1}
\\
&
-4a^{4}bKb_{1}-28ab^{2}Kb_{1}-2a^{3}b^{2}Kb_{1}+30ab^{3}Kb_{1}-18a^{2}b^{3}Kb_{1}-6aa_{12}b_{1}
\\
&
+12a^{3}a_{12}b_{1}-6ba_{12}b_{1}+48aba_{12}b_{1}-42a^{2}ba_{12}b_{1}-24a^{3}ba_{12}b_{1}-42ab^{2}a_{12}b_{1}
\\
&
+72a^{2}b^{2}a_{12}b_{1}+12b^{3}a_{12}b_{1}-24ab^{3}a_{12}b_{1}+6a^{2}a_{22}b_{1}-6a^{3}a_{22}b_{1}-18aba_{22}b_{1}
\\
&
+6a^{2}ba_{22}b_{1}+12a^{3}ba_{22}b_{1}+12b^{2}a_{22}b_{1}+18ab^{2}a_{22}b_{1}-30a^{2}b^{2}a_{22}b_{1}
\\
&
+34a^{2}b^{2}Kb_{2}-18b^{3}a_{22}b_{1}+18ab^{3}a_{22}b_{1}+6a^{2}bKb_{2}-24a^{3}bKb_{2}+6a^{4}bKb_{2}
\\
&
-4a^{4}b^{2}Kb_{2}-16ab^{3}Kb_{2}-20a^{2}b^{3}Kb_{2}+12a^{3}b^{3}Kb_{2}+8ab^{4}Kb_{2}-2a^{2}b^{4}Kb_{2}
\\
&
-6aba_{12}b_{2}-18a^{2}ba_{12}b_{2}+36a^{3}ba_{12}b_{2}+36ab^{2}a_{12}b_{2}-42a^{2}b^{2}a_{12}b_{2}
\\
&
-24a^{3}b^{2}a_{12}b_{2}-12b^{3}a_{12}b_{2}+6ab^{3}a_{12}b_{2}+30a^{2}b^{3}a_{12}b_{2}+6b^{4}a_{12}b_{2}
\\
&
-12ab^{4}a_{12}b_{2}+18a^{2}ba_{22}b_{2}-18a^{3}ba_{22}b_{2}-30ab^{2}a_{22}b_{2}+18a^{2}b^{2}a_{22}b_{2}
\\
&
+12a^{3}b^{2}a_{22}b_{2}+12b^{3}a_{22}b_{2}+6ab^{3}a_{22}b_{2}-18a^{2}b^{3}a_{22}b_{2}-6b^{4}a_{22}b_{2}%
\end{array}%
\end{equation*}

\newpage

\begin{equation*}
\begin{array}{ll}
&
+6ab^{4}a_{22}b_{2}+9a^{3}bK_{1}-3a^{4}bK_{1}-24a^{2}b^{2}K_{1}+3a^{3}b^{2}K_{1}+3a^{4}b^{2}K_{1}
\\
&
+15ab^{3}K_{1}+15a^{2}b^{3}K_{1}-12a^{3}b^{3}K_{1}-15ab^{4}K_{1}+9a^{2}b^{4}K_{1}-a^{3}bK_{2}
\\
&
+a^{4}bK_{2}+8a^{2}b^{2}K_{2}-14a^{3}b^{2}K_{2}-7ab^{3}K_{2}+13a^{2}b^{3}K_{2}+13a^{3}b^{3}K_{2}
\\
&
-a^{4}b^{3}K_{2}-20a^{2}b^{4}K_{2}+2a^{3}b^{4}K_{2}+7ab^{5}K_{2}-a^{2}b^{5}K_{2}))
\\
& -a(2(-1+b)b^{2}a_{2}^{2}(3(3a^{3}(1+b)-2b^{2}(1+b)+2ab(2+4b+b^{2}) \\
& -a^{2}(2+9b+5b^{2}))a_{11}+b(a(a^{3}(14-3b)-3b^{2}+ab(7+4b) \\
& -a^{2}(10+9b))K+3(a^{3}(-13+b)+6b^{2}-2ab(7+4b)+2a^{2}(5+9b))a_{12})) \\
& +aba_{2}(6(-1+a)(a-b)^{2}(-1+b)ba_{111}-6(-1+b)b(3a^{3}(1+b)-2b^{2}(1+b)
\\
&
+2ab(2+4b+b^{2})-a^{2}(2+9b+5b^{2}))a_{112}+24a^{2}b^{2}a_{122}-30a^{3}b^{2}a_{122}
\\
&
-42ab^{3}a_{122}+30a^{2}b^{3}a_{122}+30a^{3}b^{3}a_{122}+18b^{4}a_{122}+18ab^{4}a_{122}-54a^{2}b^{4}a_{122}
\\
&
-18b^{5}a_{122}+24ab^{5}a_{122}-2a^{3}Kb_{1}+8a^{4}Kb_{1}+12a^{2}bKb_{1}-20a^{3}bKb_{1}
\\
&
-16a^{4}bKb_{1}-4ab^{2}Kb_{1}+34a^{3}b^{2}Kb_{1}+6ab^{3}Kb_{1}-24a^{2}b^{3}Kb_{1}+6a^{3}b^{3}Kb_{1}
\\
&
-6a^{2}a_{11}b_{1}+6a^{3}a_{11}b_{1}+18aba_{11}b_{1}-6a^{2}ba_{11}b_{1}-12a^{3}ba_{11}b_{1}-12b^{2}a_{11}b_{1}
\\
&
-18ab^{2}a_{11}b_{1}+30a^{2}b^{2}a_{11}b_{1}+18b^{3}a_{11}b_{1}-18ab^{3}a_{11}b_{1}+12a^{2}a_{12}b_{1}
\\
&
-18a^{3}a_{12}b_{1}-42aba_{12}b_{1}+30a^{2}ba_{12}b_{1}+36a^{3}ba_{12}b_{1}+24b^{2}a_{12}b_{1}
\\
&
+30ab^{2}a_{12}b_{1}-84a^{2}b^{2}a_{12}b_{1}-36b^{3}a_{12}b_{1}+54ab^{3}a_{12}b_{1}-6a^{2}b^{3}a_{12}b_{1}
\\
&
-18a^{3}bKb_{2}+30a^{4}bKb_{2}-2a^{2}b^{2}Kb_{2}-28a^{4}b^{2}Kb_{2}-4ab^{3}Kb_{2}+28a^{2}b^{3}Kb_{2}
\\
&
-6a^{3}b^{3}Kb_{2}+6a^{4}b^{3}Kb_{2}+2ab^{4}Kb_{2}-14a^{2}b^{4}Kb_{2}+6a^{3}b^{4}Kb_{2}-18a^{2}ba_{11}b_{2}
\\
&
+18a^{3}ba_{11}b_{2}+30ab^{2}a_{11}b_{2}-18a^{2}b^{2}a_{11}b_{2}-12a^{3}b^{2}a_{11}b_{2}-12b^{3}a_{11}b_{2}
\\
&
-6ab^{3}a_{11}b_{2}+18a^{2}b^{3}a_{11}b_{2}+6b^{4}a_{11}b_{2}-6ab^{4}a_{11}b_{2}+48a^{2}ba_{12}b_{2}-60a^{3}ba_{12}b_{2}
\\
&
-48ab^{2}a_{12}b_{2}+30a^{2}b^{2}a_{12}b_{2}+48a^{3}b^{2}a_{12}b_{2}+24b^{3}a_{12}b_{2}-18ab^{3}a_{12}b_{2}
\\
&
-24a^{2}b^{3}a_{12}b_{2}-6a^{3}b^{3}a_{12}b_{2}-12b^{4}a_{12}b_{2}+24ab^{4}a_{12}b_{2}-6a^{2}b^{4}a_{12}b_{2}
\\
&
+a^{3}bK_{1}-7a^{4}bK_{1}-2a^{2}b^{2}K_{1}+20a^{3}b^{2}K_{1}+ab^{3}K_{1}-13a^{2}b^{3}K_{1}-13a^{3}b^{3}K_{1}
\\
&
+7a^{4}b^{3}K_{1}+14a^{2}b^{4}K_{1}-8a^{3}b^{4}K_{1}-ab^{5}K_{1}+a^{2}b^{5}K_{1}-9a^{3}b^{2}K_{2}+15a^{4}b^{2}K_{2}
\\
&
+12a^{2}b^{3}K_{2}-15a^{3}b^{3}K_{2}-15a^{4}b^{3}K_{2}-3ab^{4}K_{2}-3a^{2}b^{4}K_{2}+24a^{3}b^{4}K_{2}
\\
&
+3ab^{5}K_{2}-9a^{2}b^{5}K_{2})+6(-1+a)a((-1+b)b^{2}(-3b^{2}-a^{2}(4+b)+ab(7+b))a_{12}^{2}
\\
& +(a-b)(-1+b)b^{2}a_{11}(a(1+a-3b)K+(-1+3a-b)a_{12}-(a-b)(1+b)a_{22}) \\
& +a_{12}((-1+b)b^{3}(-a^{2}(-5+b)+a(-8+b)b+3b^{2})a_{22}+a((a+a^{2}-2b-4ab
\\
&
-2a^{2}b+6b^{2}+6ab^{2}-6b^{3})b_{1}^{2}+(a^{2}(-3+6b-4b^{2})+b^{2}(-4+8b-3b^{2})
\\
& +2ab(3-7b+4b^{2}))b_{1}b_{2}+b(b(2a^{2}+b^{2}+a(2-6b+b^{2}))b_{2}^{2} \\
& -(a-b)(-1+b)(-(1+a-3b)b_{11}+(a(3-4b)+b(-4+3b))b_{12} \\
& -b((a(7-6b)+b(-6+7b))K+(a(-3+b)+b)b_{22}))))) \\
&
+a(5a^{3}b^{2}K^{2}-10a^{2}b^{3}K^{2}-10a^{3}b^{3}K^{2}+5ab^{4}K^{2}+20a^{2}b^{4}K^{2}+5a^{3}b^{4}K^{2}
\\
&
-10ab^{5}K^{2}-10a^{2}b^{5}K^{2}+5ab^{6}K^{2}-(-1+b)b^{3}(-a^{2}(-3+b)+a(-4+b)b+b^{2})Ka_{22}
\\
&
-(a-b)^{2}(-1+b)b^{2}a_{1112}-a^{2}b^{2}a_{1122}+2ab^{3}a_{1122}-b^{4}a_{1122}+a^{2}b^{4}a_{1122}-2ab^{5}a_{1122}
\\
&
+b^{6}a_{1122}+a^{2}b^{3}a_{1222}-2ab^{4}a_{1222}-a^{2}b^{4}a_{1222}+b^{5}a_{1222}+2ab^{5}a_{1222}-b^{6}a_{1222}
\\
&
-a^{2}ba_{112}b_{1}+3ab^{2}a_{112}b_{1}+2a^{2}b^{2}a_{112}b_{1}-2b^{3}a_{112}b_{1}-5ab^{3}a_{112}b_{1}+3b^{4}a_{112}b_{1}
\\
&
+a^{2}ba_{122}b_{1}-3ab^{2}a_{122}b_{1}-2a^{2}b^{2}a_{122}b_{1}+2b^{3}a_{122}b_{1}+5ab^{3}a_{122}b_{1}-3b^{4}a_{122}b_{1}
\\
&
-a^{2}Kb_{1}^{2}-a^{3}Kb_{1}^{2}+2abKb_{1}^{2}+4a^{2}bKb_{1}^{2}+2a^{3}bKb_{1}^{2}-6ab^{2}Kb_{1}^{2}-6a^{2}b^{2}Kb_{1}^{2}
\\
&
+6ab^{3}Kb_{1}^{2}-3a^{2}b^{2}a_{112}b_{2}+5ab^{3}a_{112}b_{2}+2a^{2}b^{3}a_{112}b_{2}-2b^{4}a_{112}b_{2}-3ab^{4}a_{112}b_{2}
\\
&
+b^{5}a_{112}b_{2}+3a^{2}b^{2}a_{122}b_{2}-5ab^{3}a_{122}b_{2}-2a^{2}b^{3}a_{122}b_{2}+2b^{4}a_{122}b_{2}+3ab^{4}a_{122}b_{2}
\\
&
-b^{5}a_{122}b_{2}+3a^{3}Kb_{1}b_{2}-6a^{2}bKb_{1}b_{2}-6a^{3}bKb_{1}b_{2}+4ab^{2}Kb_{1}b_{2}+14a^{2}b^{2}Kb_{1}b_{2}
\\
&
+4a^{3}b^{2}Kb_{1}b_{2}-8ab^{3}Kb_{1}b_{2}-8a^{2}b^{3}Kb_{1}b_{2}+3ab^{4}Kb_{1}b_{2}-2a^{2}b^{2}Kb_{2}^{2}
\\
&
-2a^{3}b^{2}Kb_{2}^{2}+6a^{2}b^{3}Kb_{2}^{2}-ab^{4}Kb_{2}^{2}-a^{2}b^{4}Kb_{2}^{2}a^{2}bKb_{11}+a^{3}bKb_{11}-ab^{2}Kb_{11}%
\end{array}%
\end{equation*}

\newpage

\begin{equation*}
\begin{array}{ll}
&
-5a^{2}b^{2}Kb_{11}-a^{3}b^{2}Kb_{11}+4ab^{3}Kb_{11}+4a^{2}b^{3}Kb_{11}-3ab^{4}Kb_{11}-3a^{3}bKb_{12}
\\
&
+7a^{2}b^{2}Kb_{12}+7a^{3}b^{2}Kb_{12}-4ab^{3}Kb_{12}-14a^{2}b^{3}Kb_{12}-4a^{3}b^{3}Kb_{12}+7ab^{4}Kb_{12}
\\
&
+7a^{2}b^{4}Kb_{12}-3ab^{5}Kb_{12}-3a^{3}b^{2}Kb_{22}+4a^{2}b^{3}Kb_{22}+4a^{3}b^{3}Kb_{22}-ab^{4}Kb_{22}
\\
&
-5a^{2}b^{4}Kb_{22}-a^{3}b^{4}Kb_{22}+ab^{5}Kb_{22}+a^{2}b^{5}Kb_{22}+a^{3}bb_{1}K_{1}-3a^{2}b^{2}b_{1}K_{1}
\\
&
-2a^{3}b^{2}b_{1}K_{1}+2ab^{3}b_{1}K_{1}+5a^{2}b^{3}b_{1}K_{1}-3ab^{4}b_{1}K_{1}+3a^{3}b^{2}b_{2}K_{1}-5a^{2}b^{3}b_{2}K_{1}
\\
&
-2a^{3}b^{3}b_{2}K_{1}+2ab^{4}b_{2}K_{1}+3a^{2}b^{4}b_{2}K_{1}-ab^{5}b_{2}K_{1}-a^{3}bb_{1}K_{2}+3a^{2}b^{2}b_{1}K_{2}
\\
&
+2a^{3}b^{2}b_{1}K_{2}-2ab^{3}b_{1}K_{2}-5a^{2}b^{3}b_{1}K_{2}+3ab^{4}b_{1}K_{2}-3a^{3}b^{2}b_{2}K_{2}+5a^{2}b^{3}b_{2}K_{2}
\\
&
+2a^{3}b^{3}b_{2}K_{2}-2ab^{4}b_{2}K_{2}-3a^{2}b^{4}b_{2}K_{2}+ab^{5}b_{2}K_{2}-a^{3}b^{2}K_{11}+2a^{2}b^{3}K_{11}
\\
&
+a^{3}b^{3}K_{11}-ab^{4}K_{11}-2a^{2}b^{4}K_{11}+ab^{5}K_{11}+a^{3}b^{2}K_{12}-2a^{2}b^{3}K_{12}+ab^{4}K_{12}
\\
&
-a^{3}b^{4}K_{12}+2a^{2}b^{5}K_{12}-ab^{6}K_{12}-a^{3}b^{3}K_{22}+2a^{2}b^{4}K_{22}+a^{3}b^{4}K_{22}-ab^{5}K_{22}
\\
& -2a^{2}b^{5}K_{22}+ab^{6}K_{22})))%
\end{array}%
\end{equation*}

{\emph{Authors' addresses:} }

{\ Deparment of Mathematical Sciences, New Jersey Institute of Technology,
University Heights, Newark, NJ 07102, USA; }$\emph{E}$-$\emph{mail\;address}$%
{: vlgold@oak.njit.edu }

{\ Department of Mathematics, The University of Tromso, N9037, Tromso,
Norway; }$\emph{E}$-$\emph{mail\;address}${: lychagin@math.uit.no }

\end{document}